\documentclass[11pt]{article}
\usepackage{amssymb}
\usepackage{amsmath}
\usepackage{amsthm}
\usepackage{color}

\newtheorem{theorem}{Theorem}

\newtheorem{corollary}{Corollary}

\newtheorem{example}{Example}
\newtheorem{lemma}{Lemma}
\newtheorem{proposition}{Proposition}
\newtheorem{assumption}{Assumption}
\newtheorem{remark}{Remark}
\numberwithin{equation}{section}
\numberwithin{lemma}{section}
\numberwithin{theorem}{section}
\numberwithin{remark}{section}
\numberwithin{corollary}{section}
\numberwithin{proposition}{section}
\numberwithin{definition}{section}
\numberwithin{example}{section}
\numberwithin{assumption}{section}

\numberwithin{table}{section}

\input{tcilatex}
\begin{document}

\title{A fundamental mean-square convergence theorem for SDEs with locally
Lipschitz coefficients and its applications}
\author{M.V. Tretyakov$^{\dag }$ and Z. Zhang$^{\ddag }$ \\
$^{\dag }$School of Mathematical Sciences, University of Nottingham, \\
Nottingham, NG7 2RD, UK\\
Email: Michael.Tretyakov@nottingham.ac.uk \\
$^{\ddag }$Division of Applied Mathematics, Brown University, Providence RI,
02912\\
Email: Zhongqiang\_Zhang@brown.edu }
\maketitle

\begin{abstract}
A version of the fundamental mean-square convergence theorem is proved for
stochastic differential equations (SDE) which coefficients are allowed to
grow polynomially at infinity and which satisfy a one-sided Lipschitz
condition. The theorem is illustrated on a number of particular numerical
methods, including a special balanced scheme and fully implicit methods.
Some numerical tests are presented.


\end{abstract}

\section{Introduction}

Let $(\Omega ,\mathcal{F},P)$ be a probability space and $(w(t),\mathcal{F}%
_{t}^{w})=((w_{1}(t),\ldots ,w_{m}(t))^{\top },\mathcal{F}_{t}^{w})$ be an $%
m $-dimensional standard Wiener process, where $\mathcal{F}_{t}^{w},\ 0\leq
t\leq T,$ is an increasing family of $\sigma $-subalgebras of $\mathcal{F}$
induced by $w(t).$ We consider the system of Ito stochastic differential
equations (SDE): 
\begin{equation}
dX=a(t,X)dt+\sum_{r=1}^{m}\sigma _{r}(t,X)dw_{r}(t),\ \ t\in (t_{0},T],\
X(t_{0})=X_{0},  \label{Imps}
\end{equation}%
where $X,$\ $a,$\ $\sigma _{r}$ are $d$-dimensional column-vectors and $%
X_{0} $ is independent of $w$. We suppose that any solution $%
X_{t_{0},X_{0}}(t)$ of (\ref{Imps}) is regular on $[t_{0},T]$. We recall 
\cite{Has-B80} that a process is called regular if it is defined for all $%
t_{0}\leq t\leq T.$

In traditional numerical analysis for SDE \cite{GN,KP,MT6} it is assumed
that the SDE coefficients are globally Lipschitz which is a significant
limitation taking into account that most of the models of applicable
interest have coefficients which grow faster at infinity than a linear
function. If the global Lipschitz condition is violated, the convergence of
many usual numerical methods can disappear (see, e.g., \cite%
{Tal99,HMS,Hu,GNT04}). This has been the motivation for the recent interest
in both theoretical support of existing numerical methods and developing new
methods or approaches for solving SDE under nonglobal Lipschitz assumptions
on the coefficients.

In most of SDE applications (e.g., in molecular dynamics, financial
engineering and other problems of mathematical physics), one is interested
in simulating averages $\mathbb{E}\varphi (X(T))$ of the solution to SDE --
the task for which the weak-sense SDE approximation is sufficient and
effective \cite{GN,MT6}. The problem with divergence of weak-sense schemes
was addressed in \cite{GNT04} (see also \cite{GNTerg}) for simulation of
averages at finite time and also of ergodic limits when ensemble averaging
is used. The concept of rejecting exploding\ trajectories proposed and
justified in \cite{GNT04} allows us to use any numerical method for solving
SDE with nonglobally Lipschitz coefficients for estimating averages.
Following this concept, we do not take into account the approximate
trajectories $X(t)$ which leave a sufficiently large ball $%
S_{R}:=\{x:|x|<R\} $ during the time $T.$ See other approaches for resolving
this problem in the context of computing averages, including the case of
simulating ergodic limits via time averaging, e.g. in \cite{Tal99,Stua,Nawaf}%
.

In this paper we deal with mean-square (strong) approximation of SDE with
nonglobal Lipschitz coefficients. Mean-square schemes have their own area of
applicability (e.g. for simulating scenarios, visualization of stochastic
dynamics, filtering, etc., see further discussion on this in \cite%
{KP,MT6,HutJen12} and references therein). Furthermore, mean-square
approximation is of theoretical interest and it also provides fundamental
insight for weak-sense schemes.

We note that in the case of weak approximation we often have to simulate
large dimensional complicated stochastic systems using the Monte Carlo
technique (or time averaging), which is typical for molecular dynamics
applications, or we have to perform calculations on a daily basis, which is
usual, e.g., in financial applications. Hence the cost per step of a weak
numerical integrator should be low, which, in particular, essentially
prohibits the use of implicit methods. In contrast, areas of applicability
of mean-square schemes, as a rule, do not involve simulation of a large
number of trajectories or over very long time periods and, consequently,
there are more relaxed requirements on the cost per step of mean-square
schemes and efficient and reliable implicit schemes have practical interest.
There have been a number of recent works, including \cite%
{Hu,HMS,HutJenKlo12,SzpMao10,HutJen12} (see also the references therein),
where strong schemes for SDE with nonglobal Lipschitz coefficients were
considered. An extended literature review on this topic is available in \cite%
{HutJen12}.

In this paper we give a variant of the fundamental mean-square convergence
theorem in the case of SDE with nonglobal Lipschitz coefficients, which is
analogous to Milstein's fundamental theorem for the global Lipschitz case 
\cite{8} (see also \cite{GN,MT6}). More precisely, we assume that the SDE
coefficients can grow polynomially at infinity and satisfy a one-sided
Lipschitz condition. The theorem is stated in Section~\ref{sec:theo} and
proved in Appendix~\ref{sec:proof}. Its corollary on almost sure convergence
is also given. In Section~\ref{sec:theo} we start discussion on
applicability of the fundamental theorem, including its application to the
drift-implicit Euler scheme and thus establish its order of convergence.
Strong convergence (but without order) of this scheme was proved for SDE
with nonglobal Lipschitz drift and diffusion in \cite{SzpMao10,HutJen12}. A
particular balanced method (see the class of balanced methods in \cite%
{MPS98,MT6}) is proposed and its convergence with order $1/2$ in the
nonglobal Lipschitz setting is proved in Section~\ref{sec:tam}. In Section~%
\ref{sec:full} we revisit fully implicit (i.e., implicit both in drift and
diffusion) mean-square schemes proposed in \cite{MRT2} (see also \cite{MT6}%
). In \cite{MRT2,MT6} their convergence was proved for SDE with globally
Lipschitz coefficients. Here we relax these conditions as the drift is
required to satisfy only a one-sided Lipschitz condition and be of not
faster than polynomial growth at infinity. Some numerical experiments
supporting our results are presented in Section~\ref{sec:num}.

\section{Fundamental theorem\label{sec:theo}}

Let $X_{t_{0},X_{0}}(t)=X(t),$\ $t_{0}\leq t\leq T,$ be a solution of the
system (\ref{Imps}). In what follows we will assume the following.

\begin{assumption}
\label{asup:one-side-lip} (i) The initial condition is such that%
\begin{equation}
\mathbb{E}|X_{0}|^{2p}\leq K<\infty ,\ \ \text{for all \ }p\geq 1.
\label{mom0}
\end{equation}%
(ii) For a sufficiently large $p_{0}\geq 1$ there exists a constant $%
c_{1}\geq 0$ such that 
\begin{equation}
(x-y,a(t,x)-a(t,y))+\frac{2p_{0}-1}{2}\sum_{r=1}^{m}|\sigma _{r}(t,x)-\sigma
_{r}(t,y)|^{2}\leq c_{1}|x-y|^{2},\ t\in \lbrack t_{0},T],\ x,y\in \mathbb{R}%
^{d}.  \label{olc2}
\end{equation}%
(iii) There exist $c_{2}\geq 0$ and $\varkappa \geq 1$ such that 
\begin{equation}
|a(t,x)-a(t,y)|^{2}\leq c_{2}(1+|x|^{2\varkappa -2}+|y|^{2\varkappa
-2})|x-y|^{2},\ \ \ t\in \lbrack t_{0},T],\ x,y\in \mathbb{R}^{d}.
\label{olc3}
\end{equation}
\end{assumption}

We note that (\ref{olc2}) implies that 
\begin{equation}
(x,a(t,x))+\frac{2p_{0}-3}{2}\sum_{r=1}^{m}|\sigma _{r}(t,x)|^{2}\leq
c_{0}+c_{1}^{\prime }|x|^{2},\ \ \ t\in \lbrack t_{0},T],\ x\in \mathbb{R}%
^{d},  \label{olc}
\end{equation}%
where $c_{0}=|a(t,0)|^{2}/2+\frac{(2p_{0}-3)(2p_{0}-1)}{4}%
\sum_{r=1}^{m}|\sigma _{r}(t,0)|^{2}$ and $c_{1}^{\prime }=c_{1}+1/2.$ The
inequality (\ref{olc}) together with (\ref{mom0}) is sufficient to ensure
finiteness of moments \cite{Has-B80}: there is $K>0$ 
\begin{equation}
\mathbb{E}|X_{t_{0},X_{0}}(t)|^{2p}<K(1+\mathbb{E}|X_{0}|^{2p}),\ 1\leq
p\leq p_{0}-1,\ \ t\in \lbrack t_{0},T].  \label{Xmom}
\end{equation}%
Also, (\ref{olc3}) implies that%
\begin{equation}
|a(t,x)|^{2}\leq c_{3}+c_{2}^{\prime }|x|^{2\varkappa },\ \ \ t\in \lbrack
t_{0},T],\ x\in \mathbb{R}^{d},  \label{olc31}
\end{equation}%
where $c_{3}=2|a(t,0))|^{2}+2c_{2}(\varkappa -1)/\varkappa $ and $%
c_{2}^{\prime }=2c_{2}(1+\varkappa )/\varkappa .$

Introduce the one-step approximation $\bar{X}_{t,x}(t+h),$ $t_{0}\leq
t<t+h\leq T,$ for the solution $X_{t,x}(t+h)$ of (\ref{Imps}), which depends
on the initial point $(t,x)$, a time step\ $h,$ and $\{w_{1}(\theta
)-w_{1}(t),\ldots ,w_{m}(\theta )-w_{m}(t),\;t\leq \theta \leq t+h\}$ and
which is defined as follows: 
\begin{equation}
\bar{X}_{t,x}(t+h)=x+A(t,x,h;w_{i}(\theta )-w_{i}(t),\;i=1,\ldots ,m,\;t\leq
\theta \leq t+h).  \label{Ba03}
\end{equation}%
Using the one-step approximation (\ref{Ba03}), we recurrently construct the
approximation $(X_{k},\mathcal{F}_{t_{k}}),\;k=0,\ldots
,N,\;t_{k+1}-t_{k}=h_{k+1},\;T_{N}=T:$%
\begin{gather}
X_{0}=X(t_{0}),\;X_{k+1}=\bar{X}_{t_{k},\bar{X}_{k}}(t_{k+1})  \label{Ba04}
\\
=X_{k}+A(t_{k},X_{k},h_{k+1};w_{i}(\theta )-w_{i}(t_{k}),\;i=1,\ldots
,m,\;t_{k}\leq \theta \leq t_{k+1}).  \notag
\end{gather}

The following theorem is a generalization of Milstein's fundamental theorem 
\cite{8} (see also \cite[Chapter 1]{GN,MT6}) from the global to nonglobal
Lipschitz case. It also has similarities with a strong convergence theorem
in \cite{HMS} proved for the case of nonglobal Lipschitz drift, global
Lipschitz diffusion and Euler-type schemes.

For simplicity, we will consider a uniform time discretization, i.e. $%
h_{k}=h $ for all $k.$

\begin{theorem}
\label{thm:Bat01-lp} Suppose

(i) Assumption \ref{asup:one-side-lip} holds;

(ii) The one-step approximation $\bar{X}_{t,x}(t+h)$ from $(\ref{Ba03})$ has
the following orders of accuracy: for some $p\geq 1$ there are $\alpha \geq
1,$ $h_{0}>0,$ and $K>0$ such that for arbitrary $t_{0}\leq t\leq T-h,$\ $%
x\in \mathbf{R}^{d},$ and all $0<h\leq h_{0}:$ 
\begin{equation}
|\mathbb{E}[X_{t,x}(t+h)-\bar{X}_{t,x}(t+h)]|\leq K(1+|x|^{2\alpha
})^{1/2}h^{q_{1}}\,,  \label{Ba05-lp}
\end{equation}%
\begin{equation}
\left[ \mathbb{E}|X_{t,x}(t+h)-\bar{X}_{t,x}(t+h)|^{2p}\right] ^{1/(2p)}\leq
K(1+|x|^{2\alpha p})^{1/(2p)}h^{q_{2}}\,  \label{Ba06-lp}
\end{equation}%
with 
\begin{equation}
q_{2}\geq \frac{1}{2}\,,\;q_{1}\geq q_{2}+\frac{1}{2}\,;  \label{Ba07-lp}
\end{equation}%
(iii) The approximation $X_{k}$ from $(\ref{Ba04})$ has finite moments,
i.e., for some $p\geq 1$ there are $\beta \geq 1,$ $h_{0}>0,$ and $K>0$ such
that for all $0<h\leq h_{0}$ and all $k=0,\ldots ,N$: 
\begin{equation}
\mathbb{E}|X_{k}|^{2p}<K(1+\mathbb{E}|X_{0}|^{2p\beta }).  \label{appmomt}
\end{equation}

Then for any $N$ and $k=0,1,\ldots ,N$ the following inequality holds: 
\begin{equation}
\left[ \mathbb{E}|X_{t_{0},X_{0}}(t_{k})-\bar{X}_{t_{0},X_{0}}(t_{k})|^{2p}%
\right] ^{1/(2p)}\leq K(1+\mathbb{E}|X_{0}|^{2\gamma
p})^{1/(2p)}h^{q_{2}-1/2}\,,  \label{Ba08-lp}
\end{equation}%
where $K>0$ and $\gamma \geq 1$ do not depend on $h$ and $k,$ i.e., the
order of accuracy of the method $(\ref{Ba04})$ is $q=q_{2}-1/2.$
\end{theorem}

The theorem is proved in Appendix~\ref{sec:proof} and it uses the following
lemma.

\begin{lemma}
\label{lm:Bat02-lp} Suppose Assumption~\ref{asup:one-side-lip} holds. For
the representation 
\begin{equation}
X_{t,x}(t+\theta )-X_{t,y}(t+\theta )=x-y+Z_{t,x,y}(t+\theta ),
\label{Ba09-lp}
\end{equation}%
we have for $1\leq p\leq (p_{0}-1)/\varkappa :$ 
\begin{equation}
\mathbb{E}|X_{t,x}(t+h)-X_{t,y}(t+h)|^{2p}\leq |x-y|^{2p}(1+Kh)\,,
\label{Ba10-lp}
\end{equation}%
\begin{equation}
\mathbb{E}\left\vert Z_{t,x,y}(t+h)\right\vert ^{2p}\leq K(1+|x|^{2\varkappa
-2}+|y|^{2\varkappa -2})^{p/2}|x-y|^{2p}h^{p}\,.  \label{Ba11-lp}
\end{equation}
\end{lemma}

This lemma is proved in Appendix~\ref{sec:prooflem}. Theorem~\ref%
{thm:Bat01-lp}\ has the following corollary.

\begin{corollary}
In the setting of Theorem~\ref{thm:Bat01-lp} for $p\geq 1/(2q)$ in $(\ref%
{Ba08-lp}),$ there is $0<\varepsilon <q$ and an a.s. finite random variable $%
C(\omega )>0$ such that 
\begin{equation*}
|X_{t_{0},X_{0}}(t_{k})-X_{k}|\leq C(\omega )h^{q-\varepsilon },
\end{equation*}%
i.e., the method $(\ref{Ba04})$ for $(\ref{Imps})$ converges with order $%
q-\varepsilon $ a.s.\ 
\end{corollary}

The corollary is proved using the Borel-Cantelli-type of arguments (see,
e.g. \cite{Gyo98,filter}).

\subsection{Discussion}

In this section we make a number of observations concerning Theorem~\ref%
{thm:Bat01-lp}.

\textbf{1.} As a rule, it is not difficult to check the conditions (\ref%
{Ba05-lp})-(\ref{Ba06-lp}) following the usual routine calculations as in
the global Lipschitz case \cite{GN,KP,MT6}. We note that in order to achieve
the optimal $q_{1}$ and $q_{2}$ in (\ref{Ba05-lp})-(\ref{Ba06-lp})
additional assumptions on smoothness of $a(t,x)$ and $\sigma _{r}(t,x)$ are
usually needed.

In contrast to the conditions (\ref{Ba05-lp})-(\ref{Ba06-lp}), checking the
condition (\ref{appmomt}) on moments of a method $X_{k}$ is often rather
difficult. In the case of global Lipschitz coefficients, boundedness of
moments of $X_{k}$ is just direct implication of the boundedness of moments
of the SDE solution and the one-step properties of the method (see \cite[%
Lemma~1.1.5]{MT6}). There is no result of this type in the case of nonglobal
Lipschitz SDE and each scheme requires a special consideration. For a number
of strong schemes boundedness of moments in nonglobal Lipschitz cases were
proved (see, e.g. \cite{Hu,HMS,HutJenKlo12,HutJen12,Tal99}). In Section~\ref%
{sec:tam}\ we show boundedness of moments for a balanced method and in
Section~\ref{sec:full}\ for fully implicit methods.

Roughly speaking, Theorem~\ref{thm:Bat01-lp} says that if moments of $X_{k}$
are bounded and the scheme was proved to be convergent with order $q$ in the
global Lipschitz case then the scheme has the same convergence order $q$ in
the considered nonglobal Lipschitz case.

\textbf{2. }Assumptions and the statement of Theorem~\ref{thm:Bat01-lp}
include the famous fundamental theorem of Milstein \cite{8} proved under the
global conditions on the SDE coefficients (of course, as discussed in the
previous point, this case does not need the assumption (\ref{appmomt})).
Though the main focus here is on cases when drift and diffusion can grow
faster than a linear function at infinity, we note that the assumptions also
include the case when the diffusion coefficient grows slower than linear
function at infinity, e.g. they cover so-called CIR process which is used in
modelling short interest rates and stochastic volatility in financial
engineering.

\textbf{3.} Consider the drift-implicit scheme \cite[p. 30]{MT6}: 
\begin{equation}
X_{k+1}=X_{k}+a(t_{k+1},X_{k+1})h+\sum_{r=1}^{m}\sigma _{r}(t_{k},X_{k})\xi
_{rk}\sqrt{h},\,  \label{dimp}
\end{equation}%
where $\xi _{rk}=(w_{r}(t_{k+1})-w_{r}(t_{k}))/\sqrt{h}$ are Gaussian $%
\mathcal{N}(0,1)$ i.i.d. random variables. Assume that the coefficients $%
a(t,x)$ and $\sigma _{r}(t,x)$ have continuous first-order partial
derivatives in $t$ and the coefficient $a(t,x)$ also has continuous
first-order partial derivatives in $x^{i}$ and that all these derivatives
and the coefficients themselves satisfy inequalities of the form (\ref{olc3}%
). It is not difficult to show that the one-step approximation corresponding
to (\ref{dimp}) satisfies (\ref{Ba05-lp}) and (\ref{Ba06-lp}) with $q_{1}=2$
and $q_{2}=1,$ respectively. Its boundedness of moments, in particular,
under the condition (\ref{olc}) for time steps $h\leq 1/(2c_{1})$, is proved
in \cite{HutJen12}. Then, due to Theorem~\ref{thm:Bat01-lp}, (\ref{dimp})
converges with mean-square order $q=1/2$ (note that for $q=1/2,$ it is
sufficient to have $q_{1}=3/2$ which can be obtained under lesser smoothness
of $a).$ Further, in the case of additive noise (i.e., $\sigma
_{r}(t,x)=\sigma _{r}(t),$ $r=1,\ldots ,m),$ $q_{1}=2$ and $q_{2}=3/2$ and (%
\ref{dimp}) converges with mean-square order $1$ due to Theorem~\ref%
{thm:Bat01-lp}. We note that convergence of (\ref{dimp}) with order $1/2$ in
the global Lipschitz case is well known \cite{GN,KP,MT6}; in the case of
nonglobal Lipschitz drift and global Lipschitz diffusion was proved in \cite%
{Hu,HMS} (see also related results in \cite{Gyo98,Tal99}); and its strong
convergence without order under Assumption~\ref{asup:one-side-lip} was
proved in \cite{SzpMao10,HutJen12}.

\textbf{5.} Due to the bound (\ref{Xmom}) on the moments of the solution $%
X(t)$, it would be natural to require that $\beta $ in (\ref{appmomt})
should be equal to $1.$ Indeed, (\ref{appmomt}) with $\beta =1$ holds for
the drift-implicit method (\ref{dimp}) \cite{HutJen12} and for fully
implicit methods (see Section~\ref{sec:full}). However, this is not the case
for tamed-type methods (see \cite{HutJenKlo12}) or the balanced method from
Section~\ref{sec:tam}.

\textbf{6.} The constant $K$ in (\ref{Ba08-lp}) depends on $p,$ $t_{0},$ $T$
as well as on the SDE coefficients. The constant $\gamma $ in (\ref{Ba08-lp}%
) depends on $\alpha ,$ $\beta $ and $\varkappa .$

\section{A balanced method\label{sec:tam}}

In this section we introduce a particular balanced scheme from the class of
balanced methods introduced in \cite{MPS98} (see also \cite{MT6}) and prove
its mean-square convergence with order $1/2$ using Theorem~\ref{thm:Bat01-lp}%
. As far as we know, this variant of balanced schemes has not been
considered before. In Section~\ref{sec:num} we test the balanced scheme on a
model problem and demonstrate that it is more efficient than the tamed
scheme (\ref{fully-tamed}) (see Section~\ref{sec:num}) from \cite{HutJen12}.
We also note that it was mentioned in \cite{HutJen12} that a balanced scheme
suitable for the nonglobal Lipschitz case could potentially be derived.

Consider the following balanced-type scheme for (\ref{Imps}): 
\begin{equation}
X_{k+1}=X_{k}+\frac{a(t_{k},X_{k})h+\sum_{r=1}^{m}\sigma
_{r}(t_{k},X_{k})\xi _{rk}\sqrt{h}}{1+h|a(t_{k},X_{k})|+\sqrt{h}%
\sum_{r=1}^{m}|\sigma _{r}(t_{k},X_{k})\xi _{rk}|},  \label{ntm}
\end{equation}%
where $\xi _{rk}$ are Gaussian $\mathcal{N}(0,1)$ i.i.d. random variables.

We prove two lemmas which show that the scheme (\ref{ntm}) satisfies the
conditions of Theorem~\ref{thm:Bat01-lp}. The first lemma is on boundedness
of moments, which uses a stopping time technique (see also, e.g. \cite%
{GNT04,HutJen12}).

\begin{lemma}
\label{lem:momtm}Suppose Assumption \ref{asup:one-side-lip} holds with
sufficiently large $p_{0}$. For all natural $N$ and all $k=0,\ldots ,N$ the
following inequality holds for moments of the scheme $(\ref{ntm})$: 
\begin{equation}
\mathbb{E}|X_{k}|^{2p}\leq K(1+\mathbb{E}|X_{0}|^{2p\beta }),\ \ 1\leq p\leq 
\frac{p_{0}-1}{4(3\varkappa -2)}-\frac{1}{2},  \label{tml1}
\end{equation}%
with some constants $\beta \geq 1$ and $K>0$ independent of $h$ and $k.$
\end{lemma}

\noindent \textbf{Proof}. In the proof we shall use the letter $K$ to denote
various constants which are independent of $h$ and $k.$ We note in passing
that the case $\varkappa =1$ (i.e., when $a(t,x)$ is globally Lipschitz) is
trivial.

The following elementary consequence of the inequalities (\ref{olc}) and (%
\ref{olc31}) will be used in the proof: for any $C_{1}>0$ and $C_{2}>0:$%
\begin{eqnarray}
C_{1}\sum_{r=1}^{m}|\sigma _{r}(t,x)|^{2} &\leq &(2C_{1}c+\frac{C_{1}^{2}}{%
C_{2}})(1+|x|^{2})+C_{2}|a(t,x)|^{2}  \label{olcs} \\
&\leq &(2C_{1}c+\frac{C_{1}^{2}}{C_{2}})(1+|x|^{2})+C_{2}c\left(
1+|x|^{2\varkappa }\right) ,  \notag
\end{eqnarray}%
where $c=\max (c_{0},c_{1}^{\prime },c_{2}^{\prime },c_{3}).$

We observe that 
\begin{equation}
|X_{k+1}|\leq |X_{k}|+1\leq |X_{0}|+(k+1).  \label{tml2}
\end{equation}

Let $R>0$ be a sufficiently large number. Introduce the events 
\begin{equation}
\tilde{\Omega}_{R,k}:=\{\omega :|X_{l}|\leq R,\ l=0,\ldots ,k\},\ 
\label{event}
\end{equation}%
and their compliments $\tilde{\Lambda}_{R,k}.$ We first prove the lemma for
integer $p\geq 1.$ We have 
\begin{eqnarray}
&&\mathbb{E}\chi _{\tilde{\Omega}_{R,k+1}}(\omega )|X_{k+1}|^{2p}\leq 
\mathbb{E}\chi _{\tilde{\Omega}_{R,k}}(\omega )|X_{k+1}|^{2p}=\mathbb{E}\chi
_{\tilde{\Omega}_{R,k}}(\omega )|(X_{k+1}-X_{k})+X_{k}|^{2p}  \label{tml3} \\
&\leq &\mathbb{E}\chi _{\tilde{\Omega}_{R,k}}(\omega )|X_{k}|^{2p}+\mathbb{E}%
\chi _{\tilde{\Omega}_{R,k}}(\omega )\left\vert X_{k}\right\vert ^{2p-2}%
\left[ 2p(X_{k},X_{k+1}-X_{k})+p(2p-1)|X_{k+1}-X_{k}|^{2}\right]  \notag \\
&&+K\sum_{l=3}^{2p}\mathbb{E}\chi _{\tilde{\Omega}_{R,k}}(\omega )\left\vert
X_{k}\right\vert ^{2p-l}|X_{k+1}-X_{k}|^{l}.  \notag
\end{eqnarray}%
Consider the second term in the right-hand side of (\ref{tml3}): 
\begin{eqnarray}
&&\mathbb{E}\chi _{\tilde{\Omega}_{R,k}}(\omega )\left\vert X_{k}\right\vert
^{2p-2}\left[ 2p(X_{k},X_{k+1}-X_{k})+p(2p-1)|X_{k+1}-X_{k}|^{2}\right]
\label{tml4} \\
&=&2p\mathbb{E}\chi _{\tilde{\Omega}_{R,k}}(\omega )\left\vert
X_{k}\right\vert ^{2p-2}\mathbb{E}\left[ \left( X_{k},\frac{%
a(t_{k},X_{k})h+\sum_{r=1}^{m}\sigma _{r}(t_{k},X_{k})\xi _{rk}\sqrt{h}}{%
1+h|a(t_{k},X_{k})|+\sqrt{h}\sum_{r=1}^{m}|\sigma _{r}(t_{k},X_{k})\xi _{rk}|%
}\right) \right.  \notag \\
&&\left. \left. +\frac{2p-1}{2}\left\vert \frac{a(t_{k},X_{k})h+%
\sum_{r=1}^{m}\sigma _{r}(t_{k},X_{k})\xi _{rk}\sqrt{h}}{1+h|a(t_{k},X_{k})|+%
\sqrt{h}\sum_{r=1}^{m}|\sigma _{r}(t_{k},X_{k})\xi _{rk}|}\right\vert
^{2}\right\vert \mathcal{F}_{t_{k}}\right] .  \notag
\end{eqnarray}%
Since 
\begin{equation}
\mathbb{E}\left[ \left. \frac{\sum_{r=1}^{m}\sigma _{r}(t_{k},X_{k})\xi _{rk}%
\sqrt{h}}{1+h|a(t_{k},X_{k})|+\sqrt{h}\sum_{r=1}^{m}|\sigma
_{r}(t_{k},X_{k})\xi _{rk}|}\right\vert \mathcal{F}_{t_{k}}\right] =0
\label{tml41}
\end{equation}%
and for $l\neq r$%
\begin{equation}
\mathbb{E}\left[ \left. \frac{\sigma _{r}(t_{k},X_{k})\xi _{rk}\sqrt{h}%
\sigma _{l}(t_{k},X_{k})\xi _{lk}\sqrt{h}}{(1+h|a(t_{k},X_{k})|+\sqrt{h}%
\sum_{r=1}^{m}|\sigma _{r}(t_{k},X_{k})\xi _{rk}|)^{2}}\right\vert \mathcal{F%
}_{t_{k}}\right] =0,  \label{tml42}
\end{equation}%
the conditional expectation in (\ref{tml4}) becomes 
\begin{eqnarray}
A &:&=\mathbb{E}\left[ \left( X_{k},\frac{a(t_{k},X_{k})h+\sum_{r=1}^{m}%
\sigma _{r}(t_{k},X_{k})\xi _{rk}\sqrt{h}}{1+h|a(t_{k},X_{k})|+\sqrt{h}%
\sum_{r=1}^{m}|\sigma _{r}(t_{k},X_{k})\xi _{rk}|}\right) \right.
\label{tml5} \\
&&\left. \left. +\frac{2p-1}{2}\left\vert \frac{a(t_{k},X_{k})h+%
\sum_{r=1}^{m}\sigma _{r}(t_{k},X_{k})\xi _{rk}\sqrt{h}}{1+h|a(t_{k},X_{k})|+%
\sqrt{h}\sum_{r=1}^{m}|\sigma _{r}(t_{k},X_{k})\xi _{rk}|}\right\vert
^{2}\right\vert \mathcal{F}_{t_{k}}\right]  \notag \\
&=&\mathbb{E}\left[ \frac{\left( X_{k},a(t_{k},X_{k})h\right) }{%
1+h|a(t_{k},X_{k})|+\sqrt{h}|\sum_{r=1}^{m}\sigma _{r}(t_{k},X_{k})\xi _{rk}|%
}\right.  \notag \\
&&\left. \left. +\frac{2p-1}{2}\frac{a^{2}(t_{k},X_{k})h^{2}+h\sum_{r=1}^{m}%
\left( \sigma _{r}(t_{k},X_{k})\xi _{rk}\right) ^{2}}{\left(
1+h|a(t_{k},X_{k})|+\sqrt{h}|\sum_{r=1}^{m}\sigma _{r}(t_{k},X_{k})\xi
_{rk}|\right) ^{2}}\right\vert \mathcal{F}_{t_{k}}\right]  \notag
\end{eqnarray}%
\begin{eqnarray*}
&\leq &\mathbb{E}\left[ \frac{\left( X_{k},a(t_{k},X_{k})h\right) }{%
1+h|a(t_{k},X_{k})|+\sqrt{h}\sum_{r=1}^{m}|\sigma _{r}(t_{k},X_{k})\xi _{rk}|%
}\right. \\
&&\left. \left. +\frac{2p-1}{2}\frac{h\sum_{r=1}^{m}|\sigma
_{r}(t_{k},X_{k})|^{2}\xi _{rk}^{2}}{1+h|a(t_{k},X_{k})|+\sqrt{h}%
\sum_{r=1}^{m}|\sigma _{r}(t_{k},X_{k})\xi _{rk}|}\right\vert \mathcal{F}%
_{t_{k}}\right] +\frac{2p-1}{2}a^{2}(t_{k},X_{k})h^{2} \\
&=&\mathbb{E}\left[ \frac{\left( X_{k},a(t_{k},X_{k})h\right) +\frac{2p-1}{2}%
h\sum_{r=1}^{m}|\sigma _{r}(t_{k},X_{k})|^{2}}{1+h|a(t_{k},X_{k})|+\sqrt{h}%
\sum_{r=1}^{m}|\sigma _{r}(t_{k},X_{k})\xi _{rk}|}\right. \\
&&\left. \left. +\frac{2p-1}{2}\frac{h\sum_{r=1}^{m}|\sigma
_{r}(t_{k},X_{k})|^{2}(\xi _{rk}^{2}-1)}{1+h|a(t_{k},X_{k})|+\sqrt{h}%
\sum_{r=1}^{m}|\sigma _{r}(t_{k},X_{k})\xi _{rk}|}\right\vert \mathcal{F}%
_{t_{k}}\right] +\frac{2p-1}{2}a^{2}(t_{k},X_{k})h^{2}.
\end{eqnarray*}%
Using (\ref{olc}) and (\ref{olc31}), we obtain%
\begin{eqnarray}
A &\leq &c_{0}h+c_{1}^{\prime }|X_{k}|^{2}h\   \label{tml6} \\
&&+\frac{2p-1}{2}h\sum_{r=1}^{m}|\sigma _{r}(t_{k},X_{k})|^{2}\mathbb{E}%
\left[ \left. \frac{(\xi _{rk}^{2}-1)}{1+h|a(t_{k},X_{k})|+\sqrt{h}%
\sum_{r=1}^{m}|\sigma _{r}(t_{k},X_{k})\xi _{rk}|}\right\vert \mathcal{F}%
_{t_{k}}\right]  \notag \\
&&+Kh^{2}+K|X_{k}|^{2\varkappa }h^{2}.  \notag
\end{eqnarray}%
For the expectation in the second term in (\ref{tml6}), we obtain 
\begin{eqnarray}
&&\mathbb{E}\left[ \left. \frac{(\xi _{rk}^{2}-1)}{1+h|a(t_{k},X_{k})|+\sqrt{%
h}\sum_{r=1}^{m}|\sigma _{r}(t_{k},X_{k})\xi _{rk}|}\right\vert \mathcal{F}%
_{t_{k}}\right]  \label{tml61} \\
&=&\mathbb{E}\left[ \left. (\xi _{rk}^{2}-1)\left[ 1-\frac{h|a(t_{k},X_{k})|+%
\sqrt{h}|\sum_{l=1}^{m}\sigma _{l}(t_{k},X_{k})\xi _{lk}|}{%
1+h|a(t_{k},X_{k})|+\sqrt{h}\sum_{l=1}^{m}|\sigma _{l}(t_{k},X_{k})\xi _{lk}|%
}\right] \right\vert \mathcal{F}_{t_{k}}\right]  \notag \\
&=&-\mathbb{E}\left[ \left. (\xi _{rk}^{2}-1)\frac{h|a(t_{k},X_{k})|+\sqrt{h}%
|\sum_{r=1}^{m}\sigma _{l}(t_{k},X_{k})\xi _{lk}|}{1+h|a(t_{k},X_{k})|+\sqrt{%
h}\sum_{l=1}^{m}|\sigma _{r}(t_{k},X_{k})\xi _{lk}|}\right\vert \mathcal{F}%
_{t_{k}}\right]  \notag \\
&\leq &Kh|a(t_{k},X_{k})|+K\sqrt{h}\sum_{r=1}^{m}|\sigma _{r}(t_{k},X_{k})|.
\notag
\end{eqnarray}%
Using (\ref{olc31}) and (\ref{olcs}), we obtain from (\ref{tml6})-(\ref%
{tml61}): 
\begin{eqnarray}
A &\leq &c_{0}h+c_{1}^{\prime }|X_{k}|^{2}h+Kh\sum_{r=1}^{m}|\sigma
_{r}(t_{k},X_{k})|^{2}\left[ h|a(t_{k},X_{k})|+\sqrt{h}\sum_{r=1}^{m}|\sigma
_{r}(t_{k},X_{k})|\right]  \label{tml7} \\
&&+Kh^{2}+K|X_{k}|^{2\varkappa }h^{2}  \notag \\
&\leq &Kh(1+|X_{k}|^{2}+|X_{k}|^{2\varkappa }h+|X_{k}|^{3\varkappa
}h^{1/2})\leq Kh(1+|X_{k}|^{2}+|X_{k}|^{3\varkappa }h^{1/2}).  \notag
\end{eqnarray}

Now consider the last term in (\ref{tml3}): 
\begin{eqnarray}
&&\mathbb{E}\chi _{\tilde{\Omega}_{R,k}}(\omega )\left\vert X_{k}\right\vert
^{2p-l}|X_{k+1}-X_{k}|^{l}  \label{tml8} \\
&\leq &K\mathbb{E}\chi _{\tilde{\Omega}_{R,k}}(\omega )\left\vert
X_{k}\right\vert ^{2p-l}\left[ h^{l}|a(t_{k},X_{k})|^{l}+h^{l/2}%
\sum_{r=1}^{m}|\sigma _{r}(t_{k},X_{k})|^{l}|\xi _{rk}|^{^{l}}\right]  \notag
\\
&\leq &K\mathbb{E}\chi _{\tilde{\Omega}_{R,k}}(\omega )\left\vert
X_{k}\right\vert ^{2p-l}h^{l/2}\left[ 1+|X_{k}|^{l\varkappa }\right] , 
\notag
\end{eqnarray}%
where we used (\ref{olc31}) and (\ref{olcs}) again as well as the fact that $%
\chi _{\tilde{\Omega}_{R,k}}(\omega )$ and $X_{k}$ are $\mathcal{F}_{t_{k}}$%
-measurable while $\xi _{rk}$ are independent of $\mathcal{F}_{t_{k}}$.

Combining (\ref{tml3}), (\ref{tml4}), (\ref{tml5}), (\ref{tml7}) and (\ref%
{tml8}), we obtain 
\begin{eqnarray}
&&\mathbb{E}\chi _{\tilde{\Omega}_{R,k+1}}(\omega )|X_{k+1}|^{2p}
\label{tml9} \\
&\leq &\mathbb{E}\chi _{\tilde{\Omega}_{R,k}}(\omega )|X_{k}|^{2p}+Kh\mathbb{%
E}\chi _{\tilde{\Omega}_{R,k}}(\omega )\left\vert X_{k}\right\vert ^{2p-2}%
\left[ 1+|X_{k}|^{2}+|X_{k}|^{3\varkappa }h^{1/2}\right]  \notag \\
&&+K\sum_{l=3}^{2p}\mathbb{E}\chi _{\tilde{\Omega}_{R,k}}(\omega )\left\vert
X_{k}\right\vert ^{2p-l}h^{l/2}\left[ 1+|X_{k}|^{l\varkappa }\right]  \notag
\\
&\leq &\mathbb{E}\chi _{\tilde{\Omega}_{R,k}}(\omega )|X_{k}|^{2p}+Kh\mathbb{%
E}\chi _{\tilde{\Omega}_{R,k}}(\omega )\left\vert X_{k}\right\vert
^{2p}+K\sum_{l=2}^{2p}\mathbb{E}\chi _{\tilde{\Omega}_{R,k}}(\omega
)\left\vert X_{k}\right\vert ^{2p-l}h^{l/2}  \notag \\
&&+Kh^{3/2}\mathbb{E}\chi _{\tilde{\Omega}_{R,k}}(\omega )\left\vert
X_{k}\right\vert ^{2p-2+3\varkappa }+Kh\sum_{l=3}^{2p}\mathbb{E}\chi _{%
\tilde{\Omega}_{R,k}}(\omega )\left\vert X_{k}\right\vert ^{2p+l(\varkappa
-1)}h^{l/2-1}.  \notag
\end{eqnarray}%
Choosing 
\begin{equation}
R=R(h)=h^{-1/(6\varkappa -4)},  \label{tml10}
\end{equation}%
we get $\mathbb{E}\chi _{\tilde{\Omega}_{R,k}}(\omega )\left\vert
X_{k}\right\vert ^{2p-2+3\varkappa }h^{l/2-1}\leq \chi _{\tilde{\Omega}%
_{R(h),k}}(\omega )\left\vert X_{k}\right\vert ^{2p}$ and $\chi _{\tilde{%
\Omega}_{R(h),k}}(\omega )\left\vert X_{k}\right\vert ^{2p+l(\varkappa
-1)}h^{l/2-1}\leq \chi _{\tilde{\Omega}_{R(h),k}}(\omega )\left\vert
X_{k}\right\vert ^{2p},$ $l=3,\ldots ,2p,$ and hence we re-write (\ref{tml9}%
) as%
\begin{eqnarray}
&&\mathbb{E}\chi _{\tilde{\Omega}_{R(h),k+1}}(\omega )|X_{k+1}|^{2p}
\label{tml11} \\
&\leq &\mathbb{E}\chi _{\tilde{\Omega}_{R(h),k}}(\omega )|X_{k}|^{2p}+Kh%
\mathbb{E}\chi _{\tilde{\Omega}_{R(h),k}}(\omega )\left\vert
X_{k}\right\vert ^{2p}+K\sum_{l=1}^{p}\mathbb{E}\chi _{\tilde{\Omega}%
_{R(h),k}}(\omega )\left\vert X_{k}\right\vert ^{2(p-l)}h^{l}  \notag \\
&\leq &\mathbb{E}\chi _{\tilde{\Omega}_{R(h),k}}(\omega )|X_{k}|^{2p}+Kh%
\mathbb{E}\chi _{\tilde{\Omega}_{R(h),k}}(\omega )\left\vert
X_{k}\right\vert ^{2p}+Kh,  \notag
\end{eqnarray}%
where in the last line we have used Young's inequality. From here, we get by
Gronwall's inequality that 
\begin{equation}
\mathbb{E}\chi _{\tilde{\Omega}_{R(h),k}}(\omega )|X_{k}|^{2p}\leq K(1+%
\mathbb{E}|X_{0}|^{2p}),  \label{tml12}
\end{equation}%
where $R(h)$ is from (\ref{tml10}) and $K$ does not depend on $k$ and $h$
but it depends on $p.$

It remains to estimate $\mathbb{E}\chi _{\tilde{\Lambda}_{R(h),k}}(\omega
)|X_{k}|^{2p}.$ We have 
\begin{eqnarray*}
\chi _{\tilde{\Lambda}_{R,k}} &=&1-\chi _{\tilde{\Omega}_{R,k}}=1-\chi _{%
\tilde{\Omega}_{R,k-1}}\chi _{|X_{k}|\leq R}=\chi _{\tilde{\Lambda}%
_{R,k-1}}+\chi _{\tilde{\Omega}_{R,k-1}}\chi _{|X_{k}|>R} \\
&=&\cdots =\sum_{l=0}^{k}\chi _{\tilde{\Omega}_{R,l-1}}\chi _{|X_{l}|>R},
\end{eqnarray*}%
where we put $\chi _{\tilde{\Omega}_{R,-1}}=1.$ Then, using (\ref{tml2}), (%
\ref{tml12}), (\ref{mom0}), and Cauchy-Bunyakovsky's and Markov's
inequalities, we obtain%
\begin{gather}
\mathbb{E}\chi _{\tilde{\Lambda}_{R(h),k}}(\omega )|X_{k}|^{2p}=\mathbb{E}%
\sum_{l=0}^{k}|X_{k}|^{2p}\chi _{\tilde{\Omega}_{R(h),l-1}}\chi
_{|X_{l}|>R(h)}  \label{Markov} \\
\leq \left( \mathbb{E}|X_{0}+k|^{4p}\right) ^{1/2}\sum_{l=0}^{k}\left( 
\mathbb{E}\left[ \chi _{\tilde{\Omega}_{R(h),l-1}|X_{l}|>R(h)}\right]
\right) ^{1/2}  \notag \\
=\left( \mathbb{E}|X_{0}+k|^{4p}\right) ^{1/2}\sum_{l=0}^{k}\left( P(\chi _{%
\tilde{\Omega}_{R(h),l-1}}|X_{l}|>R)\right) ^{1/2}  \notag \\
\leq \left( \mathbb{E}|X_{0}+k|^{4p}\right) ^{1/2}\sum_{l=0}^{k}\frac{\left( 
\mathbb{E}(\chi _{\tilde{\Omega}_{R(h),l-1}}|X_{l}|^{2(2p+1)(6\varkappa
-4)})\right) ^{1/2}}{R(h)^{(2p+1)(6\varkappa -4)}}  \notag \\
\leq K\left( \mathbb{E}|X_{0}+k|^{4p}\right) ^{1/2}\left( \mathbb{E}%
(1+|X_{0}|^{2(2p+1)(6\varkappa -4)})\right) ^{1/2}kh^{2p+1}\leq K(1+\mathbb{E%
}|X_{0}|^{4p+2(2p+1)(6\varkappa -4)})^{1/2},  \notag
\end{gather}%
which together with (\ref{tml12}) implies (\ref{tml1}) for integer $p\geq 1$%
. Then, by Jensen's inequality, (\ref{tml1}) holds for non-integer $p$ as
well. $\square $

The next lemma gives estimates for the one-step error of the balanced scheme
(\ref{ntm}).

\begin{lemma}
\label{lem:loctm}Assume that $(\ref{Xmom})$ holds. Assume that the
coefficients $a(t,x)$ and $\sigma _{r}(t,x)$ have continuous first-order
partial derivatives in $t$ and that these derivatives and the coefficients
satisfy inequalities of the form $(\ref{olc3})$. Then the scheme $(\ref{ntm}%
) $ satisfies the inequalities $(\ref{Ba05-lp})$ and $(\ref{Ba06-lp})$ with $%
q_{1}=3/2$ and $q_{2}=1,$ respectively.
\end{lemma}

The proof of this lemma is a routine analysis of the one-step approximation
corresponding to (\ref{ntm}) using the equalities (\ref{tml41})-(\ref{tml42}%
). Since such analysis is similar to those done in the global Lipschitz case 
\cite{MPS98,MT6}, we omit these routine calculations here. Lemmas~\ref%
{lem:momtm} and~\ref{lem:loctm} and Theorem~\ref{thm:Bat01-lp} imply the
following result.

\begin{proposition}
\label{prp:tm}Under the assumptions of Lemmas~\ref{lem:momtm} and~\ref%
{lem:loctm} the balanced scheme $(\ref{ntm})$ has mean-square order $1/2,$
i.e., for it the inequality $(\ref{Ba08-lp})$ holds with $q=q_{2}-1/2=1/2.$
\end{proposition}

\begin{remark}
In the additive noise case the mean-square order of the balanced scheme $(%
\ref{ntm})$ does not improve ($q_{1}$ and $q_{2}$ remain $3/2$ and $1,$
respectively).
\end{remark}

\section{Fully implicit schemes\label{sec:full}}

Fully implicit (i.e., implicit both in drift and diffusion coefficients)
mean-square schemes were proposed in \cite{MRT2} (see also \cite[Chapter 1]%
{MT6}), where their convergence was proved under global Lipschitz
conditions. Here we analyze these schemes under the following assumptions,
which are stronger with respect to the diffusion coefficient than Assumption~%
\ref{asup:one-side-lip}\ used in the previous Sections~\ref{sec:theo}\ and~%
\ref{sec:tam}.

\begin{assumption}
\label{asup:nonglAglobS} (i) The initial condition is such that%
\begin{equation}
\mathbb{E}|X_{0}|^{2p}\leq K<\infty ,\ \ \text{for all \ }p\geq 1.
\label{as20}
\end{equation}%
(ii) There exists a constant $c_{1}\geq 0$ such that 
\begin{equation}
(x-y,a(t,x)-a(t,y))\leq c_{1}|x-y|^{2},\ \ t\in \lbrack t_{0},T],\ x,y\in 
\mathbb{R}^{d}.  \label{as21}
\end{equation}%
(iii) There exist $c_{2}\geq 0$ and $\varkappa \geq 1$ such that 
\begin{equation}
|a(t,x)-a(t,y)|^{2}\leq c_{2}(1+|x|^{2\varkappa -2}+|y|^{2\varkappa
-2})|x-y|^{2},\ \ t\in \lbrack t_{0},T],\ x,y\in \mathbb{R}^{d}.\ 
\label{as22}
\end{equation}%
(iv) The coefficients $\sigma _{r}(t,x)$ have continuous bounded first-order
spatial derivatives so that there are constants $L_{1}\geq 0$ and $L_{2}\geq
0:$ 
\begin{equation}
|\nabla \sigma _{r}(t,x)|\leq L_{1},\ r=1,\ldots ,m,\ \ t\in \lbrack
t_{0},T],\ x\in \mathbb{R}^{d},  \label{as23}
\end{equation}%
and 
\begin{equation}
|\nabla \sigma _{r}(t,x)\sigma _{r}(t,x)-\nabla \sigma _{r}(t,y)\sigma
_{r}(t,y)|\leq L_{2}|x-y|,\ ,\ r=1,\ldots ,m,\ \ t\in \lbrack t_{0},T],\
x,y\in \mathbb{R}^{d}.\   \label{as233}
\end{equation}
\end{assumption}

In proofs which follow we will need some implications of Assumption~\ref%
{asup:nonglAglobS}. The condition (\ref{as21}) implies that there is $c\geq
0 $%
\begin{equation}
(x,a(t,x))\leq c(1+|x|^{2}),\ \ \ t\in \lbrack t_{0},T],\ x\in \mathbb{R}%
^{d}.  \label{as24}
\end{equation}%
It follows from (\ref{as23}) that 
\begin{equation}
|\sigma _{r}(t,x)-\sigma _{r}(t,y)|\leq L_{1}|x-y|,\ \ t\in \lbrack
t_{0},T],\ x,y\in \mathbb{R}^{d},  \label{as255}
\end{equation}%
and hence 
\begin{equation}
|\sigma _{r}(t,x)|\leq L_{1}|x|+L_{0},  \label{as257}
\end{equation}%
where $L_{0}=\max_{t\in \lbrack t_{0},T]}|\sigma _{r}(t,0)|.$ Further, there
is $L\geq 0:$%
\begin{equation}
|\nabla \sigma _{r}(t,x)\sigma _{r}(t,x)|\leq L(1+|x|),\ \ \ t\in \lbrack
t_{0},T],\ x\in \mathbb{R}^{d},  \label{as259}
\end{equation}%
and 
\begin{equation}
|\sigma _{r}(t,x)|^{2}\leq L(1+|x|^{2}),\ \ \ t\in \lbrack t_{0},T],\ x\in 
\mathbb{R}^{d}.  \label{as25}
\end{equation}

For definiteness, we consider the following one-parametric family of methods
for (\ref{Imps}) from the broader class of fully implicit schemes of \cite%
{MRT2,MT6}: 
\begin{eqnarray}
X_{k+1} &=&X_{k}+a(t_{k+\lambda },(1-\lambda )X_{k}+\lambda X_{k+1})h\ 
\label{IEM} \\
&&-\lambda \sum_{r=1}^{m}\sum_{j=1}^{d}\frac{\partial \sigma _{r}}{\partial
x^{j}}(t_{k+\lambda },(1-\lambda )X_{k}+\lambda X_{k+1})\sigma
_{r}^{j}(t_{k+\lambda },(1-\lambda )X_{k}+\lambda X_{k+1})h  \notag \\
&&+\sum_{r=1}^{m}\sigma _{r}(t_{k+\lambda },(1-\lambda )X_{k}+\lambda
X_{k+1})\left( \zeta _{rh}\right) _{k}\sqrt{h},  \notag
\end{eqnarray}%
where $0\leq \lambda \leq 1,$ $t_{k+\lambda }=t_{k}+\lambda h$ and $\left(
\zeta _{rh}\right) _{k}$ are i.i.d. random variables so that 
\begin{equation}
\zeta _{h}=\left\{ 
\begin{array}{c}
\xi ,\;|\xi |\leq A_{h}, \\ 
A_{h},\;\xi >A_{h}, \\ 
-A_{h},\;\xi <-A_{h},%
\end{array}%
\right.  \label{Fin61}
\end{equation}%
with $\xi $ $\sim $ $\mathcal{N}(0,1)$ and $A_{h}=\sqrt{2l|\ln h|}$ with $%
l\geq 1.$ We recall \cite[Lemma 1.3.4]{MT6} that 
\begin{equation}
E(\xi ^{2}-\zeta _{h}^{2})=(1+2\sqrt{2l|\ln h|})h^{l}.  \label{fu1}
\end{equation}

\begin{remark}
Three choices of $\lambda $ are most notable: $\lambda =0$ gives the
explicit Euler scheme which is divergent \cite{HMS,HutJen12} in the
considered setting; $\lambda =1$ gives the fully implicit Euler scheme; and $%
\lambda =1/2$ corresponds to the mid-point rule, which in application to a
system of Stratonovich SDE is derivative free \cite[p. 45]{MT6}.
\end{remark}

Now we will study properties of the method (\ref{IEM}).

Consider the one-step approximations corresponding to (\ref{IEM}) 
\begin{eqnarray}
\bar{X} &=&\bar{X}^{\lambda }=x+a(t+\lambda h,U^{\lambda })h-\lambda
\sum_{r=1}^{m}\sum_{j=1}^{d}\frac{\partial \sigma _{r}}{\partial x^{j}}%
(t+\lambda h,U^{\lambda })\sigma _{r}^{j}(t+\lambda h,U^{\lambda })h\ 
\label{IEM1} \\
&&+\sum_{r=1}^{m}\sigma _{r}(t+\lambda h,U^{\lambda })\zeta _{rh}\sqrt{h}, 
\notag
\end{eqnarray}%
where 
\begin{equation}
U=U^{\lambda }:=(1-\lambda )x+\lambda \bar{X}^{\lambda }.  \label{U}
\end{equation}%
Note that 
\begin{eqnarray}
U^{\lambda } &=&x+\lambda a(t+\lambda h,U^{\lambda })-\lambda
^{2}\sum_{r=1}^{m}\sum_{j=1}^{d}\frac{\partial \sigma _{r}}{\partial x^{j}}%
(t+\lambda h,U^{\lambda })\sigma _{r}^{j}(t+\lambda h,U^{\lambda })h\ 
\label{MPM2} \\
&&+\lambda \sum_{r=1}^{m}\sigma _{r}(t+\lambda h,U^{\lambda })\zeta _{rh}%
\sqrt{h}.  \notag
\end{eqnarray}

\begin{lemma}
\label{lem:res}Let $0<\lambda \leq 1.$ Assume that Assumption~\ref%
{asup:nonglAglobS} holds. For an arbitrary $0<\varepsilon <1,$ find $h_{0}>0$
such that 
\begin{equation}
\lambda \left[ h_{0}c_{1}+m\lambda L_{2}h_{0}+mL_{1}\sqrt{2lh_{0}|\ln h_{0}|}%
\right] =1-\varepsilon .  \label{IEMh0}
\end{equation}%
Then the equation $(\ref{IEM1})$ for any $0<h\leq h_{0}$ has the unique
solution $\bar{X}$\ which satisfies the inequalities for some $K>0:$ 
\begin{equation}
|\bar{X}-x|\leq K(1+|x|^{\varkappa })h+K(1+|x|)\sqrt{h|\ln h|}  \label{Fim16}
\end{equation}%
and%
\begin{equation}
|\bar{X}|^{2}\leq \frac{16}{3\varepsilon ^{2}\lambda }(L_{0}+1)\sqrt{2lh|\ln
h|}+\frac{4}{\lambda ^{2}}[(1-\lambda )^{2}+\frac{4}{3\varepsilon ^{2}}%
]|x|^{2},\ t\in \lbrack t_{0},T],\ x\in \mathbb{R}^{d}.\   \label{lres0}
\end{equation}
\end{lemma}

\noindent \textbf{Proof}. Let 
\begin{equation}
\tilde{a}(t,x)=a(t,x)-\lambda \sum_{r=1}^{m}\sum_{j=1}^{d}\frac{\partial
\sigma _{r}}{\partial x^{j}}(t,x)\sigma _{r}^{j}(t,x).  \label{atild}
\end{equation}%
For any fixed $\lambda ,$ $t,$\ $\zeta _{rh},$ and $h,$ we introduce the
function 
\begin{equation*}
\psi (z)=z-\lambda \tilde{a}(t+\lambda h,z)h-\lambda \sum_{r=1}^{m}\sigma
_{r}(t+\lambda h,z)\zeta _{rh}\sqrt{h}
\end{equation*}%
which is continuous in $z$ due to our assumptions. The equation (\ref{IEM1})
can be written as 
\begin{equation}
\psi (U^{\lambda })=x.  \label{lres2}
\end{equation}%
Using (\ref{as21}), (\ref{as233}) and (\ref{as255}), we obtain 
\begin{eqnarray}
(z-y,\psi (z)-\psi (y)) &\geq &|z-y|^{2}-h\lambda c_{1}|z-y|^{2}-hm\lambda
^{2}L_{2}|z-y|^{2}\ \   \label{lres3} \\
&&-m\lambda L_{1}\sqrt{2lh|\ln h|}|z-y|^{2}  \notag \\
&=&(1-\lambda \left[ hc_{1}+m\lambda L_{2}h+mL_{1}\sqrt{2lh|\ln h|}\right]
)|z-y|^{2}\geq \varepsilon |z-y|^{2}>0,  \notag
\end{eqnarray}%
i.e., $\psi (z)$ is uniformly monotone function for $h\leq h_{0}$. This
implies (see, e.g. \cite[Theorem 6.4.4, p. 167]{Orteg}) that (\ref{IEM1})
has a unique solution.

We obtain from (\ref{lres2}) and (\ref{lres3}): 
\begin{eqnarray*}
\varepsilon |U|^{2} &\leq &(U,\psi (U)-\psi (0))=(U,x-\psi (0)) \\
&\leq &\frac{\varepsilon }{4}|U|^{2}+\frac{2}{\varepsilon }|x|^{2}+\frac{2}{%
\varepsilon }|\psi (0)|^{2}\leq \frac{\varepsilon }{4}|U|^{2}+\frac{2}{%
\varepsilon }|x|^{2}+\frac{2\lambda (L_{0}+1)\sqrt{2lh|\ln h|}}{\varepsilon }%
,
\end{eqnarray*}%
from which (\ref{lres0}) follows.%
\begin{equation}
|U|^{2}\leq \frac{8}{3\varepsilon ^{2}}(\lambda (L_{0}+1)\sqrt{2lh|\ln h|}%
+|x|^{2})  \label{Ubound}
\end{equation}%
which implies (\ref{lres0}).

Further, it follows from (\ref{U}), (\ref{lres2}) and (\ref{lres3}) that%
\begin{eqnarray*}
\lambda \varepsilon |\bar{X}-x|^{2} &=&\varepsilon |U-x|^{2}\leq
(U-x,-\lambda \tilde{a}(t+\lambda h,x)h-\lambda \sum_{r=1}^{m}\sigma
_{r}(t+\lambda h,x)\zeta _{rh}\sqrt{h}) \\
&\leq &\lambda ^{2}|\bar{X}-x|\left( h|\tilde{a}(t+\lambda h,x)|+\sqrt{%
2lh|\ln h|}\sum_{r=1}^{m}|\sigma _{r}(t+\lambda h,x)|\right) .
\end{eqnarray*}%
Then, using (\ref{as22}) and (\ref{as257}), we obtain (\ref{Fim16}), which
completes the proof of Lemma~\ref{lem:res} for the implicit method (\ref{IEM}%
). $\ \square $

Now we consider boundedness of moments of (\ref{IEM}).

\begin{lemma}
\label{lem:momfie}Let $1/2<\lambda \leq 1.$ Assume that Assumption~\ref%
{asup:nonglAglobS} holds. Then for all $0<h\leq h_{0}$ with $h_{0}$ from $(%
\ref{IEMh0})$ and for all $k=0,\ldots ,N$ the following inequality holds for
the fully implicit scheme $(\ref{IEM})$ for $p\geq 1$: 
\begin{equation}
\mathbb{E}|X_{k}|^{2p}\leq K(1+\mathbb{E}|X_{0}|^{2p}),  \label{liem1}
\end{equation}%
where $K>0$ is a constant.
\end{lemma}

\noindent \textbf{Proof}. We note that (\ref{as24}) and (\ref{as259}) imply%
\begin{equation}
(x,\tilde{a}(t,x))\leq (c+3m\lambda L/2)(1+|x|^{2}),\ \ t\in \lbrack
t_{0},T],\ x\in \mathbb{R},  \label{as26}
\end{equation}%
which together with (\ref{as23}) ensures that the solution of (\ref{Imps})
has all moments (\ref{Xmom}), $p\geq 1$ \cite{Has-B80}.

Let $U_{k+1}=(1-\lambda )X_{k}+\lambda X_{k+1}$ (cf. (\ref{U})). We have 
\begin{gather}
V_{k+1}:=|X_{k+1}|^{2}-|X_{k}|^{2}=2(U_{k+1},X_{k+1}-X_{k})-(2\lambda
-1)|X_{k+1}-X_{k}|^{2}  \label{liem2} \\
=2\lambda h(U_{k+1},\tilde{a}(t_{k+\lambda },U_{k+1}))+2\lambda \sqrt{h}%
(U_{k+1},\sum_{r=1}^{m}\sigma _{r}(t_{k+\lambda },U_{k+1})\left( \zeta
_{rh}\right) _{k})  \notag \\
-(2\lambda -1)h^{2}|\tilde{a}(t_{k+\lambda },U_{k+1})|^{2}-(2\lambda
-1)h|\sum_{r=1}^{m}\sigma _{r}(t_{k+\lambda },U_{k+1})\left( \zeta
_{rh}\right) _{k}|^{2}  \notag \\
-2(2\lambda -1)h^{3/2}(\tilde{a}(t_{k+\lambda
},U_{k+1}),\sum_{r=1}^{m}\sigma _{r}(t_{k+\lambda },U_{k+1})\left( \zeta
_{rh}\right) _{k})  \notag \\
=2\lambda (U_{k+1},h\tilde{a}(t_{k+\lambda },U_{k+1}))+2\lambda \sqrt{h}%
(X_{k},\sum_{r=1}^{m}\sigma _{r}(t_{k+\lambda },U_{k+1})\left( \zeta
_{rh}\right) _{k})  \notag \\
+(2\lambda ^{2}-2\lambda +1)h|\sum_{r=1}^{m}\sigma _{r}(t_{k+\lambda
},U_{k+1})\left( \zeta _{rh}\right) _{k}|^{2}-(2\lambda -1)h^{2}|\tilde{a}%
(t_{k+\lambda },U_{k+1})|^{2}  \notag \\
+2(1-\lambda )^{2}h^{3/2}(\tilde{a}(t_{k+\lambda
},U_{k+1}),\sum_{r=1}^{m}\sigma _{r}(t_{k+\lambda },U_{k+1})\left( \zeta
_{rh}\right) _{k}).  \notag
\end{gather}%
Expanding $\sigma _{r}(t_{k+\lambda },U_{k+1})$ at $(t_{k+\lambda },X_{k}),$
we obtain 
\begin{gather}
V_{k+1}=2\lambda (U_{k+1},h\tilde{a}(t_{k+\lambda },U_{k+1}))+2\lambda \sqrt{%
h}(X_{k},\sum_{r=1}^{m}\sigma _{r}(t_{k+\lambda },X_{k})\left( \zeta
_{rh}\right) _{k})  \label{liem3} \\
+2\lambda \sqrt{h}(X_{k},\sum_{r=1}^{m}\nabla \sigma _{r}(t_{k+\lambda
},\theta )(U_{k+1}-X_{k})\left( \zeta _{rh}\right) _{k})  \notag \\
+(2\lambda ^{2}-2\lambda +1)h|\sum_{r=1}^{m}\sigma _{r}(t_{k+\lambda
},U_{k+1})\left( \zeta _{rh}\right) _{k}|^{2}-(2\lambda -1)h^{2}|\tilde{a}%
(t_{k+\lambda },U_{k+1})|^{2}  \notag \\
+2(\lambda -1)^{2}h^{3/2}(\tilde{a}(t_{k+\lambda
},U_{k+1}),\sum_{r=1}^{m}\sigma _{r}(t_{k+\lambda },U_{k+1})\left( \zeta
_{rh}\right) _{k})  \notag
\end{gather}%
\begin{gather*}
=2\lambda (U_{k+1},h\tilde{a}(t_{k+\lambda },U_{k+1}))+2\lambda \sqrt{h}%
(X_{k},\sum_{r=1}^{m}\sigma _{r}(t_{k+\lambda },X_{k})\left( \zeta
_{rh}\right) _{k}) \\
+2\lambda ^{2}h^{3/2}(X_{k},\sum_{r=1}^{m}\nabla \sigma _{r}(t_{k+\lambda
},\theta )\tilde{a}(t_{k+\lambda },U_{k+1})\left( \zeta _{rh}\right) _{k}) \\
+2\lambda ^{2}h(X_{k},\sum_{r=1}^{m}\nabla \sigma _{r}(t_{k+\lambda },\theta
)\sum_{l=1}^{m}\sigma _{l}(t_{k+\lambda },X_{k})\left( \zeta _{lh}\right)
_{k}\left( \zeta _{rh}\right) _{k}) \\
+(2\lambda ^{2}-2\lambda +1)h|\sum_{r=1}^{m}\sigma _{r}(t_{k+\lambda
},U_{k+1})\left( \zeta _{rh}\right) _{k}|^{2}-(2\lambda -1)h^{2}|\tilde{a}%
(t_{k+\lambda },U_{k+1})|^{2} \\
+2(1-\lambda )^{2}h^{3/2}(\tilde{a}(t_{k+\lambda
},U_{k+1}),\sum_{r=1}^{m}\sigma _{r}(t_{k+\lambda },U_{k+1})\left( \zeta
_{rh}\right) _{k}),
\end{gather*}%
where $\theta =\nu U_{k+1}-(1-\nu )X_{k},$ $\nu \in \lbrack 0,1],$ is an
intermediate point. Using (\ref{as26}), (\ref{as25}), Young's inequality, (%
\ref{as23}) and (\ref{lres0}), we obtain 
\begin{eqnarray}
V_{k+1} &\leq &\lambda h(2c+3\lambda mL)(1+|U_{k+1}|^{2})+2\lambda \sqrt{h}%
(X_{k},\sum_{r=1}^{m}\sigma _{r}(t_{k+\lambda },X_{k})\left( \zeta
_{rh}\right) _{k})  \label{blow} \\
&&+\frac{2\lambda -1}{2}h^{2}|\tilde{a}(t_{k+\lambda },U_{k+1})|^{2}+\frac{%
2\lambda ^{4}}{2\lambda -1}h|X_{k}|^{2}m\sum_{r=1}^{m}|\nabla \sigma
_{r}(t_{k+\lambda },\theta )|^{2}|\left( \zeta _{rh}\right) _{k}|^{2}  \notag
\\
&&+\lambda ^{2}hm|X_{k}|^{2}\sum_{r=1}^{m}|\nabla \sigma _{r}(t_{k+\lambda
},\theta )|^{2}|\left( \zeta _{rh}\right) _{k}|^{2}+\lambda
^{2}hm\sum_{l=1}^{m}|\sigma _{l}(t_{k+\lambda },X_{k})|^{2}|\left( \zeta
_{lh}\right) _{k}|^{2}  \notag \\
&&+(2\lambda ^{2}-2\lambda +1)hm\sum_{r=1}^{m}|\sigma _{r}(t_{k+\lambda
},U_{k+1})|^{2}|\left( \zeta _{rh}\right) _{k}|^{2}-(2\lambda -1)h^{2}|%
\tilde{a}(t_{k+\lambda },U_{k+1})|^{2}  \notag \\
&&+\frac{2\lambda -1}{2}h^{2}|\tilde{a}(t_{k+\lambda },U_{k+1})|^{2}+\frac{%
2(1-\lambda )^{4}}{2\lambda -1}hm\sum_{r=1}^{m}|\sigma _{r}(t_{k+\lambda
},U_{k+1})|^{2}|\left( \zeta _{rh}\right) _{k}|^{2}  \notag \\
&\leq &\lambda h(2c+3\lambda mL)(1+|U_{k+1}|^{2})+2\lambda \sqrt{h}%
(X_{k},\sum_{r=1}^{m}\sigma _{r}(t_{k+\lambda },X_{k})\left( \zeta
_{rh}\right) _{k})  \notag \\
&&+\lambda ^{2}[\frac{2\lambda ^{2}}{2\lambda -1}+1]hL_{1}^{2}|X_{k}|^{2}m%
\sum_{r=1}^{m}|\left( \zeta _{rh}\right) _{k}|^{2}+\lambda
^{2}hmL(1+|X_{k}|^{2})\sum_{l=1}^{m}|\left( \zeta _{lh}\right) _{k}|^{2} 
\notag \\
&&+[2\lambda ^{2}-2\lambda +1+\frac{2(1-\lambda )^{4}}{2\lambda -1}%
]hmL(1+|U_{k+1}|^{2})\sum_{l=1}^{m}|\left( \zeta _{lh}\right) _{k}|^{2}. 
\notag
\end{eqnarray}%
Then using (\ref{Ubound}), we arrive at 
\begin{equation*}
V_{k+1}\leq Kh(1+|X_{k}|^{2})\left( 1+\sum_{r=1}^{m}|\left( \zeta
_{rh}\right) _{k}|^{2}\right) +2\lambda \sqrt{h}(X_{k},\sum_{r=1}^{m}\sigma
_{r}(t_{k+\lambda },X_{k})\left( \zeta _{rh}\right) _{k}),
\end{equation*}%
where $K>0$ is independent of $h$ and $k$ while it depends on $\lambda $ and
on constants appearing in (\ref{as21})-(\ref{as25}).

Thus 
\begin{equation*}
1+|X_{k+1}|^{2}\leq 1+|X_{k}|^{2}+Kh(1+|X_{k}|^{2})\left(
1+\sum_{r=1}^{m}|\left( \zeta _{rh}\right) _{k}|^{2}\right) +2\lambda \sqrt{h%
}(X_{k},\sum_{r=1}^{m}\sigma _{r}(t_{k+\lambda },X_{k})\left( \zeta
_{rh}\right) _{k}).
\end{equation*}%
Then for integer $p\geq 1$ we get 
\begin{eqnarray*}
\left( 1+|X_{k+1}|^{2}\right) ^{p} &\leq &\left( 1+|X_{k}|^{2}\right)
^{p}+K\left( 1+|X_{k}|^{2}\right) ^{p}\sum_{l=1}^{p}h^{l}\left[
1+\sum_{r=1}^{m}|\left( \zeta _{rh}\right) _{k}|^{2}\right] ^{l} \\
&&+K\sum_{l=1}^{p}\left( 1+|X_{k}|^{2}\right) ^{p-l}h^{l/2}\left[
(X_{k},\sum_{r=1}^{m}\sigma _{r}(t_{k+\lambda },X_{k})\left( \zeta
_{rh}\right) _{k})\right] ^{l},
\end{eqnarray*}%
whence, observing that $X_{k}$ are $\mathcal{F}_{t_{k}}$-measurable while $%
\xi _{rk}$ are independent of $\mathcal{F}_{t_{k}},$ it is not difficult to
obtain 
\begin{eqnarray*}
\mathbb{E}\left( 1+|X_{k+1}|^{2}\right) ^{p} &\leq &\mathbb{E}\left(
1+|X_{k}|^{2}\right) ^{p}+Kh\mathbb{E}\left( 1+|X_{k}|^{2}\right) ^{p} \\
&&+K\sum_{l=2}^{p}\mathbb{E}\left( 1+|X_{k}|^{2}\right) ^{p-l}h^{l/2}\left[
(X_{k},\sum_{r=1}^{m}\sigma _{r}(t_{k+\lambda },X_{k})\left( \zeta
_{rh}\right) _{k})\right] ^{l} \\
&\leq &\mathbb{E}\left( 1+|X_{k}|^{2}\right) ^{p}+Kh\mathbb{E}\left(
1+|X_{k}|^{2}\right) ^{p},
\end{eqnarray*}%
which together with Gronwall's inequality completes the proof of the lemma
for integer $p\geq 1$.\ Then by Jensen's inequality for non-integer $p>1$ as
well. $\square $

We have not succeeded in proving boundedness of moments for the mid-point
scheme, i.e., $(\ref{IEM})$ with $\lambda =1/2$ under Assumption~\ref%
{asup:nonglAglobS}. One can observe that the proof of Lemma~\ref{lem:momfie}
is not applicable to this choice of $\lambda $ as the estimate in (\ref{blow}%
) blows up when $\lambda \rightarrow 1/2$ and it is clear that the mid-point
scheme is the boundary case. We also know \cite{HW} that for $\sigma _{r}=0$
(\ref{IEM}) is B-stable for $\lambda \geq 1/2$ and not B-stable (in fact,
not A-stable) for $\lambda <1/2.$ It is natural to expect that for $\lambda
<1/2$ the moments of (\ref{IEM}) are not bounded and hence the method with $%
\lambda <1/2$ is divergent under Assumption~\ref{asup:nonglAglobS} (see also
such a conclusion for the drift-implicit $\theta $-method in \cite{SzpMao10}%
). In our experiments (Section~\ref{sec:num}) the mid-point method produced
accurate results.

At the same time, we proved boundedness of moments for the mid-point scheme
if in addition to Assumption~\ref{asup:nonglAglobS} we require that the
diffusion coefficients $\sigma _{r}(t,x)$ are bounded. The proof is similar
to the proof of Lemma~\ref{lem:momtm}.

\begin{lemma}
\label{lem:momMPM}Let the assumptions of Lemma~\ref{lem:momfie} hold and in
addition assume that the diffusion coefficients $\sigma _{r}(t,x)$ are
uniformly bounded. Then the moments of the mid-point method $(\ref{IEM})$
with $\lambda =1/2$ has bounded moments: for $p\geq 1$: 
\begin{equation}
\mathbb{E}|X_{k}|^{2p}\leq K(1+\mathbb{E}|X_{0}|^{4(p+1)\varkappa -4})^{1/2},
\label{mM0}
\end{equation}%
where $K>0$ is a constant.
\end{lemma}

\noindent \textbf{Proof.} For $\varkappa =1$ (cf. (\ref{as22})), i.e., the
global Lipschitz case, boundedness of moments of $X_{k}$ is established in 
\cite{MRT2,MT6}. Let $\varkappa >1.$

From (\ref{liem2}), we have 
\begin{equation}
|X_{k+1}|^{2}-|X_{k}|^{2}=h(U_{k+1},\tilde{a}(t_{k+1/2},U_{k+1}))+\sqrt{h}%
(U_{k+1},\sum_{r=1}^{m}\sigma _{r}(t_{k+1/2},U_{k+1})\left( \zeta
_{rh}\right) _{k}).  \notag
\end{equation}%
Then using (\ref{as26}) and Young's inequality, boundedness of $\sigma _{r}$
and (\ref{Ubound}), we get%
\begin{eqnarray*}
|X_{k+1}|^{2}-|X_{k}|^{2} &\leq &h(c+3\lambda
mL/2)(1+|U_{k+1}|^{2})+h|U_{k+1}|^{2}+K\sum_{r=1}^{m}|\left( \zeta
_{rh}\right) _{k}|^{2} \\
&\leq &Kh(1+|X_{k}|^{2})+K\sum_{r=1}^{m}|\left( \zeta _{rh}\right) _{k}|^{2},
\end{eqnarray*}%
from which one can obtain for integer $p\geq 1:$ 
\begin{equation}
\mathbb{E}\left( 1+|X_{k}|^{2}\right) ^{p}\leq Kh^{-p}\mathbb{E}\left(
1+|X_{0}|^{2}\right) ^{p}.  \label{mM1}
\end{equation}%
Further, using (\ref{as26}), Young's inequality, boundedness of $\sigma
_{r}, $ (\ref{as23}), (\ref{as22}) and (\ref{Ubound}), we get from (\ref%
{liem3}): 
\begin{gather*}
|X_{k+1}|^{2}-|X_{k}|^{2}=(U_{k+1},h\tilde{a}(t_{k+1/2},U_{k+1}))+\sqrt{h}%
(X_{k},\sum_{r=1}^{m}\sigma _{r}(t_{k+1/2},X_{k})\left( \zeta _{rh}\right)
_{k}) \\
+\frac{1}{2}h^{3/2}(X_{k},\sum_{r=1}^{m}\nabla \sigma _{r}(t_{k+1},\theta )%
\tilde{a}(t_{k+1/2},U_{k+1})\left( \zeta _{rh}\right) _{k}) \\
+\frac{1}{2}h(X_{k},\sum_{r=1}^{m}\nabla \sigma _{r}(t_{k+1/2},\theta
)\sum_{l=1}^{m}\sigma _{l}(t_{k+1/2},X_{k})\left( \zeta _{lh}\right)
_{k}\left( \zeta _{rh}\right) _{k})+\frac{1}{2}h|\sum_{r=1}^{m}\sigma
_{r}(t_{k+1/2},U_{k+1})\left( \zeta _{rh}\right) _{k}|^{2} \\
+\frac{1}{2}h^{3/2}(\tilde{a}(t_{k+\lambda },U_{k+1}),\sum_{r=1}^{m}\sigma
_{r}(t_{k+\lambda },U_{k+1})\left( \zeta _{rh}\right) _{k})
\end{gather*}%
\begin{eqnarray*}
&\leq &Kh(1+|X_{k}|^{2})(1+\sum_{r=1}^{m}|\left( \zeta _{rh}\right)
_{k}|^{2})+\sqrt{h}(X_{k},\sum_{r=1}^{m}\sigma _{r}(t_{k+1/2},X_{k})\left(
\zeta _{rh}\right) _{k}) \\
&&+Kh^{3/2}(1+|X_{k}|^{\varkappa +1})\sum_{r=1}^{m}|\left( \zeta
_{rh}\right) _{k}|.
\end{eqnarray*}%
Choosing $R(h)=h^{-1/2(\varkappa -1)},$ we get after some additional
calculation (cf. (\ref{tml12})):%
\begin{equation}
\mathbb{E}\chi _{\tilde{\Omega}_{R(h),k}}(\omega )|X_{k}|^{2p}\leq K(1+%
\mathbb{E}|X_{0}|^{2p}),  \label{mM2}
\end{equation}%
where $\tilde{\Omega}_{R(h),k}$ is the event as in (\ref{event}).

Now using (\ref{mM1}), (\ref{mM2}), and Cauchy-Bunyakovsky's and Markov's
inequalities, we arrive at (cf. (\ref{Markov})): 
\begin{equation*}
\mathbb{E}\chi _{\tilde{\Lambda}_{R(h),k}}(\omega )|X_{k}|^{2p}\leq K(1+%
\mathbb{E}|X_{0}|^{4(p+1)\varkappa -4})^{1/2},
\end{equation*}%
from which together with (\ref{mM2}) the inequality (\ref{mM0}) follows. \ $%
\square $

The next lemma gives estimates for the one-step error of the method (\ref%
{IEM}).

\begin{lemma}
\label{lem:locIEM}Let $0\leq \lambda \leq 1.$ Assume that $(\ref{Xmom})$
holds. Assume that the coefficient $a(t,x)$ has continuous first order
partial derivative in $t$ and in $x^{i}$ and that the derivatives and the
coefficient satisfy inequalities of the form $(\ref{as22})$; the functions $%
\sigma _{r}(t,x)$ have continuous first-order partial derivatives in $t$ and
that the derivatives and the coefficients satisfy inequalities of the form $(%
\ref{as23})-(\ref{as233});$ and the functions $\nabla \sigma _{r}(t,x)\sigma
_{r}(t,x)$ have continuous first partial derivatives in $t$ and in $x^{i}$
which satisfy inequalities of the form $(\ref{as233}).$ Then the method $(%
\ref{IEM})$ satisfies the inequalities $(\ref{Ba05-lp})$ and $(\ref{Ba06-lp}%
) $ with $q_{1}=2$ and $q_{2}=1,$ respectively.
\end{lemma}

Proofs of this lemma is rather routine and similar to the global Lipschitz
case \cite{MRT2,MT6} and it is omitted here. Using Lemmas~\ref{lem:res}-\ref%
{lem:locIEM}, the next proposition follows from Theorem~\ref{thm:Bat01-lp}.

\begin{proposition}
\label{prp:imp}Let for $1/2<\lambda \leq 1$ the assumptions of Lemmas~\ref%
{lem:momfie} and~\ref{lem:locIEM} hold and for $\lambda =1/2$ in addition
assume that the diffusion coefficients $\sigma _{r}(t,x)$ are uniformly
bounded. Then the fully implicit method $(\ref{IEM})$ has mean-square order $%
1/2,$ i.e., for it the inequality $(\ref{Ba08-lp})$ holds with $q=1/2.$
\end{proposition}

\begin{remark}
\label{Rem:commu}Consider the commutative case, i.e., when $\Lambda
_{i}\sigma _{r}=\Lambda _{r}\sigma _{i}$ (here the operator $\Lambda
_{r}:=(\sigma _{r},\partial /\partial x))$ or in the case of a system with
one noise (i.e., $m=1).$ Then in the setting of Lemma~\ref{lem:locIEM}, the
mid-point method, i.e., $(\ref{IEM})$ with $\lambda =1/2,$ satisfies the
inequalities $(\ref{Ba05-lp})$ and $(\ref{Ba06-lp})$ with $q_{1}=2$ and $%
q_{2}=3/2,$ respectively (see such a result in the global Lipschitz case in 
\cite{MRT2,MT6}). Therefore, it converges in this case with mean-square
order $1$ when its moments are bounded.
\end{remark}

\section{Numerical examples\label{sec:num}}

In this section we will test the following schemes: the balanced method (\ref%
{ntm}) from Section~\ref{sec:tam}; the drift-implicit scheme (\ref{dimp});
the fully implicit Euler scheme (\ref{IEM}) with $\lambda =1;$ the mid-point
method (\ref{IEM}) with $\lambda =1/2;$ the drift-tamed Euler scheme (a
modified balanced method) \cite{HutJenKlo12}: 
\begin{equation}
X_{k+1}=X_{k}+h\frac{a(X_{k})}{1+h\left\vert a(X_{k})\right\vert }%
+\sum_{r=1}^{m}\sigma _{r}(t_{k},X_{k})\xi _{rk}\sqrt{h};  \label{drift_tam}
\end{equation}%
the fully-tamed scheme \cite{HutJen12}: 
\begin{equation}
X_{k+1}=X_{k}+\frac{a(X_{k})h+\sum_{r=1}^{m}\sigma _{r}(t_{k},X_{k})\xi _{rk}%
\sqrt{h}}{\max \left( {1,h\left\vert ha(X_{k})+\sum_{r=1}^{m}\sigma
_{r}(t_{k},X_{k})\xi _{rk}\sqrt{h}\right\vert }\right) };\ 
\label{fully-tamed}
\end{equation}%
and the trapezoidal scheme: 
\begin{equation}
X_{k+1}=X_{k}+\frac{h}{2}\left[ a(X_{k+1})+a(X_{k})\right]
+\sum_{r=1}^{m}\sigma _{r}(t_{k},X_{k})\xi _{rk}\sqrt{h}.  \label{trap}
\end{equation}%
As before, $\xi _{rk}=(w_{r}(t_{k+1})-w_{r}(t_{k}))/\sqrt{h}$ are Gaussian $%
\mathcal{N}(0,1)$ i.i.d. random variables. We note that under Assumption~\ref%
{asup:one-side-lip} boundedness of second moments and strong convergence
(without giving order) of $\theta $-schemes, and in particular of (\ref{trap}%
), can be found in \cite{SzpMao10}. Strong convergence with order $1/2$ of (%
\ref{drift_tam}) under Assumption~\ref{asup:nonglAglobS} is proved in \cite%
{HutJenKlo12}. Strong convergence of (\ref{fully-tamed}) without order under
Assumption~\ref{asup:one-side-lip} is proved in \cite{HutJen12}.

In all the experiments with fully implicit schemes, where the truncated
random variables $\zeta $ are used, we took $l=2$ (see (\ref{fu1})). The
experiments were performed using Matlab R2012a on a Macintosh desktop
computer with Intel Xeon CPU E5462 (quad-core, $2.80$~GHz). In simulations
we used the Mersenne twister random generator with seed $100$. Newton's
method was used to solve the nonlinear algebraic equations at each step of
the implicit schemes.

We test the methods on the two model problems. The first one satisfies
Assumption~\ref{asup:nonglAglobS} (nonglobal Lipschitz drift, global
Lipschitz diffusion) and has two non-commutative noises. The second example
satisfies Assumption~\ref{asup:one-side-lip} (nonglobal Lipschitz both drift
and diffusion). The aim of the tests is to compare performance of the
methods: their accuracy (i.e., roughly speaking, size of prefactors at a
power of $h)$ and computational costs. We note that experiments cannot prove
or disprove boundedness of moments of the schemes since experiments rely on
a finite sample of trajectories run over a finite time interval while
blow-up of moments in divergent methods (e.g., explicit Euler scheme) is, in
general, a result of large deviations \cite{Stua,GNT04}.

\begin{example}
\label{exa1}Our first test model is the Stratonovich SDE of the form: 
\begin{equation}
dX=(1-X^{5})\,dt+X\circ \,dw_{1}+dw_{2},\quad X(0)=0.  \label{Example1}
\end{equation}
\end{example}

\noindent In Ito's sense, the drift of the equation becomes $%
a(t,x)=1-x^{5}+x/2$. Here we tested the balanced method (\ref{ntm}); the
drift-tamed scheme (\ref{drift_tam}); the fully implicit Euler scheme (\ref%
{IEM}) with $\lambda =1;$ the mid-point method (\ref{IEM}) with $\lambda
=1/2.$ We note that for all the methods tested on this example except the
mid-point rule mean-square convergence with order $1/2$ is proved either in
earlier papers or here as it was described before.

To compute the mean-square error, we run $M$ independent trajectories $%
X^{(i)}(t),$ $X_{k}^{(i)}$: 
\begin{equation}
\left( E\left[ X(T)-X_{N}\right] ^{2}\right) ^{1/2}\doteq \left( \frac{1}{M}%
\sum_{i=1}^{M}[X^{(i)}(T)-X_{N}^{(i)}]^{2}\right) ^{1/2}.  \label{experr}
\end{equation}%
We took time $T=50$ and $M=10^{4}.$ The reference solution was computed by
the mid-point method with small time step $h=10^{-4}.$ It was verified that
using a different implicit scheme for simulating a reference solution does
not affect the outcome of the tests. We chose the mid-point scheme as a
reference since in all the experiments it produced the most accurate results.

Table~\ref{tab:Example1} gives the mean-square errors and experimentally
observed convergence rates for the corresponding methods. We checked that
the number of trajectories $M=10^{4}$ was sufficiently large for the
statistical errors not to significantly hinder the mean-square errors (the
Monte Carlo error computed with $95\%$ confidence was at least 10 time
smaller than the reported mean-square errors except values for (\ref%
{drift_tam}) at $h=0.1$ and $0.05$ where it was at least $5$ time smaller
than the mean-square errors). In addition to the data in the table, we
evaluated errors for (\ref{ntm}) for smaller time steps: $h=0.002$ -- the
error is $9.27$e-02 (rate $0.41$), $0.001$ -- $6.86$e-02 ($0.44$). The
observed rates of convergence of all the tested methods are close to the
predicted $1/2.$ For a fixed time step $h,$ the most accurate scheme is the
mid-point one, the less accurate scheme is the new balanced method (\ref{ntm}%
). To produce the result with accuracy\textbf{\ }$\sim $\textbf{\ }$%
0.06-0.07,$\textbf{\ }in our experiment of running $M=10^{4}$ trajectories%
\textbf{\ \ }the scheme (\ref{drift_tam}) required 170 sec.,\textbf{\ }the
mid-point (\ref{IEM}) with $\lambda =1/2$ -- 329 sec., (\ref{IEM}) with $%
\lambda =1$\ -- 723 sec., and\textbf{\ }(\ref{ntm}) -- 1870 sec.\textbf{\ }%
That is, our experiments confirmed the conclusion of \cite{HutJenKlo12} that
the drift-tamed (modified balance method) (\ref{drift_tam}) from \cite%
{HutJenKlo12} is highly competitive. We note that (\ref{drift_tam}) is not
applicable when diffusion grows faster than a linear function and that in
this case the balanced method (\ref{ntm}) can outcompete implicit schemes as
it is shown in the next example.

\begin{table}[h] \centering%
\caption{{\it Example \ref{exa1}.}
Mean-square errors of the selected schemes. See further details in the text.
\label{tab:Example1}}%
\setlength{\tabcolsep}{2pt}\centering 
\begin{tabular}{c|cc|cc|cc|cc}
\hline\hline
$h$ & (\ref{IEM}), $\lambda =1$ & rate & (\ref{IEM}), $\lambda =1/2$ & rate
& (\ref{drift_tam}) & rate & (\ref{ntm}) & rate \\ \hline
0.1 & 1.712e-01 & -- & 1.443e-01 & -- & 3.748e-01 & -- & 3.594e-01 & -- \\ 
\hline
0.05 & 1.234e-01 & 0.47 & 9.224e-02 & 0.65 & 2.103e-01 & 0.83 & 3.017e-01 & 
0.25 \\ \hline
0.02 & 7.692e-02 & 0.52 & 5.261e-02 & 0.61 & 9.472e-02 & 0.87 & 2.297e-01 & 
0.30 \\ \hline
0.01 & 5.478e-02 & 0.49 & 3.549e-02 & 0.57 & 6.104e-02 & 0.63 & 1.778e-01 & 
0.37 \\ \hline
0.005 & 3.935e-02 & 0.48 & 2.487e-02 & 0.51 & 3.959e-02 & 0.62 & 1.354e-01 & 
0.39 \\ \hline\hline
\end{tabular}
\end{table}%

\begin{example}
\label{exa2}Consider the SDE in the Stratonovich sense: 
\begin{equation}
dX=(1-X^{5})\,dt+X^{2}\circ \,dw,\quad X(0)=0.  \label{Example2}
\end{equation}
\end{example}

\noindent In Ito's sense, the drift of the equation becomes $%
a(t,x)=1-x^{5}+x^{3}$.

Here we tested the balanced method (\ref{ntm}); the fully-tamed Euler scheme
(\ref{fully-tamed}); the drift-implicit scheme (\ref{dimp}); the fully
implicit Euler scheme (\ref{IEM}) with $\lambda =1;$ the mid-point method (%
\ref{IEM}) with $\lambda =1/2;$ and the trapezoidal scheme (\ref{trap}). We
recall that in the case of nonglobal Lipschitz drift and diffusion, for the
drift-implicit scheme (\ref{dimp}) and the balanced method (\ref{ntm})
mean-square convergence with order $1/2$ is shown earlier in this paper;
strong convergence of the trapezoidal scheme (\ref{trap}) without order is
proved in \cite{SzpMao10}, it is natural to expect that its mean-square
order is $1/2$ which is indeed supported by the experiments. Strong
convergence of (\ref{fully-tamed}) without order is proved in \cite{HutJen12}%
. We note that it can be proved directly that implicit algebraic equations
arising from application of the mid-point and fully implicit Euler schemes
to (\ref{Example2}) have unique solutions under a sufficiently small time
step.

The reference solution was computed by the mid-point method with small time
step $h=10^{-4}.$ The time $T=50$ and $M=10^{4}$ in (\ref{experr}).

The fully-tamed scheme (\ref{fully-tamed}) did not produce accurate results
until the time step size is at least $h=0.005$ and we do not then report its
errors here but see the remark below.

\begin{remark}
\label{RemO}The fully-tamed scheme $(\ref{fully-tamed})$ appears to be of a
low practical value. If at a step $k_{\ast },$ the event $O:={\left\vert
ha(X_{k})+\sum_{r=1}^{m}\sigma _{r}(t_{k},X_{k})\xi _{rk}\sqrt{h}\right\vert
>1/h}$ happens, then in the case of $(\ref{Example2})$ the trajectory $%
X_{k}, $ $k>k_{\ast },$ oscillates approximately between $X_{k_{\ast }}$ and 
$X_{k_{\ast }}-signum(X_{k_{\ast }})/h$. Since the probability of the event $%
O $ is positive for any step size $h>0$ and grows with integration time, it
is unavoidable that in some scenarios (i.e., on some trajectories) such
oscillatory behavior will appear. For instance, in this experiment for $%
h=0.1 $ we observed $989$ out of $1000$ paths for which $O$ happened over
the time interval $[0,50];$ for $h=0.05$ -- $866$ out of $1000$ paths. From
the practical point of view, $(\ref{fully-tamed})$ works as long as the
explicit Euler scheme works (cf. \cite{Stua} and also \cite[p. 17]{MT6}).
The strong convergence (without order) of $(\ref{fully-tamed})$ \cite%
{HutJen12} in comparison with the explicit Euler scheme is due to the
following fact. When event $O$ happens for the Euler scheme its sequence $%
X_{k}$ starts oscillating with growing amplitude which leads to
unboundedness of its moments and, consequently, its divergence in the
mean-square sense. For $(\ref{fully-tamed})$, the oscillations are bounded
by $\sim 1/h$ and since the probability of $O$ over a finite time interval
decreasing with decrease of $h,$ then the moments are bounded uniformly in $%
h $. At the same time, the one-step approximation of $(\ref{fully-tamed})$
does not satisfy the conditions $(\ref{Ba05-lp})$ and $(\ref{Ba06-lp})$ of
Theorem~\ref{thm:Bat01-lp}. We note that the explicit balanced-type scheme $(%
\ref{ntm})$ does not have such drawbacks as $(\ref{fully-tamed})$.
\end{remark}

Table~\ref{tbl:ms-compare-errors} gives the mean-square errors and
experimentally observed convergence rates for the corresponding methods. We
checked that the number of trajectories $M=10^{4}$ was sufficiently large
for the statistical errors not to significantly hinder the mean-square
errors (the Monte Carlo error computed with $95\%$ confidence was at least
ten time smaller than the reported mean-square errors). In addition to the
data in the table, we evaluated errors for (\ref{ntm}) for smaller time
steps: $h=0.002$ -- the error is $3.70$e-02 (rate $0.41$), $0.001$ -- $2.73$%
e-02 ($0.44$), $0.0005$ -- $2.00$e-02 ($0.45$), i.e., for smaller $h$ the
observed convergence rate of (\ref{ntm}) becomes closer to the theoretically
predicted order $1/2.$ Since (\ref{Example2}) is with single noise, Remark~%
\ref{Rem:commu}\ is valid here which explains why the mid-point scheme
demonstrates the first order of convergence. The other implicit schemes show
the order $1/2$ as expected. Table~\ref{tbl:ms-compare-timing} presents the
time costs in seconds. Let us fix tolerance level at $0.05-0.06.$ We
highlighted in bold the corresponding values in both tables. We see that in
this example the mid-point scheme is the most efficient which is due to its
first order convergence in the commutative case. Among methods of order $%
1/2, $ the balanced method (\ref{ntm}) is the fastest and one can expect
that for multi-dimensional SDE the explicit scheme (\ref{ntm}) can
considerably outperform implicit methods (see a similar outcome for the
drift-tamed method (\ref{drift_tam}) supported by experiments in \cite%
{HutJenKlo12}; note that (\ref{drift_tam}), in comparison with (\ref{ntm}),
is, as a rule, divergent when diffusion is growing faster than a linear
function on infinity).

\begin{table}[h] \centering%
\caption{{\it Example \ref{exa2}.}
Mean-square errors of the selected schemes. See further details in the text.
\label{tbl:ms-compare-errors}}%
\setlength{\tabcolsep}{2pt}\centering 
\begin{tabular}{c|cc|cc|cc|cc|cc}
\hline\hline
$h$ & (\ref{dimp}) & rate & (\ref{IEM}), $\lambda =1$ & rate & (\ref{IEM}), $%
\lambda =1/2$ & rate & (\ref{trap}) & rate & (\ref{ntm}) & rate \\ \hline
0.2 & 3.449e-01 & -- & 1.816e-01 & -- & 1.378e-01 & -- & 4.920e-01 & -- & 
2.102e-01 & -- \\ \hline
0.1 & 2.441e-01 & 0.50 & 1.331e-01 & 0.45 & 8.723e-02 & 0.66 & 3.526e-01 & 
0.48 & 1.637e-01 & 0.36 \\ \hline
0.05 & 1.592e-01 & 0.62 & 9.619e-02 & 0.47 & \textbf{5.344e-02} & 0.71 & 
2.230e-01 & 0.66 & 1.270e-01 & 0.37 \\ \hline
0.02 & 8.360e-02 & 0.70 & 6.599e-02 & 0.41 & 2.242e-02 & 0.95 & 1.048e-01 & 
0.82 & 9.170e-02 & 0.36 \\ \hline
0.01 & \textbf{5.460e-02} & 0.61 & \textbf{4.919e-02} & 0.42 & 1.145e-02 & 
0.97 & \textbf{5.990e-02} & 0.81 & 7.065e-02 & 0.38 \\ \hline
0.005 & 3.682e-02 & 0.57 & 3.522e-02 & 0.48 & 5.945e-03 & 0.95 & 3.784e-02 & 
0.66 & \textbf{5.393e-02} & 0.39 \\ \hline\hline
\end{tabular}
\end{table}%

\begin{table}[h] \centering%
\caption{{\it Example \ref{exa2}.}
Comparison of computational times for the selected schemes. See further details in the text.
\label{tbl:ms-compare-timing}}%
\centering%
\begin{tabular}{c|ccccc}
\hline\hline
$h$ & (\ref{dimp}) & (\ref{IEM}), $\lambda =1$ & (\ref{IEM}), $\lambda =1/2$
& (\ref{trap}) & (\ref{ntm}) \\ \hline
0.2 & 9.25e+00 & 1.10e+01 & 9.33e+00 & 1.20e+01 & 3.98e+00 \\ \hline
0.1 & 1.77e+01 & 2.17e+01 & 1.80e+01 & 2.30e+01 & 7.49e+00 \\ \hline
0.05 & 3.42e+01 & 4.26e+01 & \textbf{3.51e+01} & 4.48e+01 & 1.41e+01 \\ 
\hline
0.02 & 8.33e+01 & 1.04e+02 & 8.69e+01 & 1.10e+02 & 3.37e+01 \\ \hline
0.01 & \textbf{1.64e+02} & \textbf{2.05e+02} & 1.73e+02 & \textbf{2.19e+02}
& 6.62e+01 \\ \hline
0.005 & 3.25e+02 & 4.07e+02 & 3.47e+02 & 4.37e+02 & \textbf{1.32e+02} \\ 
\hline\hline
\end{tabular}%
\end{table}%

\section*{Acknowledgments}

We are very grateful to Boris Rozovskii and George Karniadakis for fruitful
discussions. MVT was partially supported by the Leverhulme Trust Fellowship
SAF-2012-006 and is also grateful to ICERM (Brown University, Providence)
for its hospitality. ZZ was supported by the OSD/MURI grant
FA9550-09-1-0613.

\appendix{}

\section{Proof of the fundamental theorem \label{sec:proof}}

Note that in this and the next section we shall use the letter $K$ to denote
various constants which are independent of $h$ and $k.$ The proof exploits
the idea of the prove of this theorem in the global Lipschitz case \cite{8}.

Consider the error of the method $\bar{X}_{t_{0},X_{0}}(t_{k+1})$ at the $%
(k+1)$-step:%
\begin{gather}
\rho _{k+1}:=X_{t_{0},X_{0}}(t_{k+1})-\bar{X}%
_{t_{0},X_{0}}(t_{k+1})=X_{t_{k},X(t_{k})}(t_{k+1})-\bar{X}%
_{t_{k},X_{k}}(t_{k+1})  \label{Ba27-lp} \\
=(X_{t_{k},X(t_{k})}(t_{k+1})-X_{t_{k},X_{k}}(t_{k+1}))+(X_{t_{k},X_{k}}(t_{k+1})-%
\bar{X}_{t_{k},X_{k}}(t_{k+1}))\,.  \notag
\end{gather}%
The first difference in the right-hand side of (\ref{Ba27-lp}) is the error
of the solution arising due to the error in the initial data at time $t_{k},$
accumulated at the $k$-th step, which we can re-write as 
\begin{eqnarray*}
S_{t_{k},X(t_{k}),X_{k}}(t_{k+1})
&=&S_{k+1}:=X_{t_{k},X(t_{k})}(t_{k+1})-X_{t_{k},X_{k}}(t_{k+1})=\rho
_{k}+Z_{t_{k},X(t_{k}),X_{k}}(t_{k+1}) \\
&=&\rho _{k}+Z_{k+1},
\end{eqnarray*}%
where $Z$ is as in (\ref{Ba09-lp}). The second difference in (\ref{Ba27-lp})
is the one-step error at the $(k+1)$-step and we denote it as $r_{k+1}:$%
\begin{equation*}
r_{k+1}=X_{t_{k},X_{k}}(t_{k+1})-\bar{X}_{t_{k},X_{k}}(t_{k+1}).
\end{equation*}%
Let $p\geq 1$ be an integer. We have 
\begin{eqnarray}
\mathbb{E}|\rho _{k+1}|^{2p} &=&\mathbb{E}\left\vert
S_{k+1}+r_{k+1}\right\vert ^{2p}=\mathbb{E}%
[(S_{k+1},S_{k+1})+2(S_{k+1},r_{k+1})+(r_{k+1},r_{k+1})]^{p}\ \ \ \ \ \ 
\label{thm1} \\
&\leq &\mathbb{E}\left\vert S_{k+1}\right\vert ^{2p}+2p\mathbb{E}\left\vert
S_{k+1}\right\vert ^{2p-2}(\rho _{k}+Z_{k+1},r_{k+1})+K\sum_{l=2}^{2p}%
\mathbb{E}\left\vert S_{k+1}\right\vert ^{2p-l}|r_{k+1}|^{l}.  \notag
\end{eqnarray}%
Due to (\ref{Ba10-lp}) of Lemma~\ref{lm:Bat02-lp}, the first term on the
right-hand side of (\ref{thm1}) is estimated as%
\begin{equation}
\mathbb{E}\left\vert S_{k+1}\right\vert ^{2p}\leq \mathbb{E}|\rho
_{k}|^{2p}(1+Kh).  \label{thm2}
\end{equation}

Consider the second term on the right-hand side of (\ref{thm1}): 
\begin{gather}
\mathbb{E}\left\vert S_{k+1}\right\vert ^{2p-2}(\rho _{k}+Z_{k+1},r_{k+1})=%
\mathbb{E}\left\vert \rho _{k}\right\vert ^{2p-2}(\rho _{k},r_{k+1})
\label{thm3} \\
+\mathbb{E}\left( \left\vert S_{k+1}\right\vert ^{2p-2}-\left\vert \rho
_{k}\right\vert ^{2p-2}\right) (\rho _{k},r_{k+1})+\mathbb{E}\left\vert
S_{k+1}\right\vert ^{2p-2}(Z_{k+1},r_{k+1}).  \notag
\end{gather}%
Due to $\mathcal{F}_{t_{k}}$-measurability of $\rho _{k}$ and due to the
conditional variant of (\ref{Ba05-lp}), we get for the first term on the
right-hand side of (\ref{thm3}):%
\begin{equation}
\mathbb{E}\left\vert \rho _{k}\right\vert ^{2p-2}(\rho _{k},r_{k+1})\leq K%
\mathbb{E}\left\vert \rho _{k}\right\vert ^{2p-1}(1+|X_{k}|^{2\alpha
})^{1/2}h^{q_{1}}.  \label{thm4}
\end{equation}%
Consider the second term on the right-hand side of (\ref{thm3}) and first of
all note that it is equal to zero for $p=1$. We have for integer $p\geq 2:$%
\begin{equation*}
\mathbb{E}\left( \left\vert S_{k+1}\right\vert ^{2p-2}-\left\vert \rho
_{k}\right\vert ^{2p-2}\right) (\rho _{k},r_{k+1})\leq K\mathbb{E}\left\vert
Z_{k+1}\right\vert |\rho
_{k}||r_{k+1}|\sum_{l=0}^{2p-3}|S_{k+1}|^{2p-3-l}|\rho _{k}|^{l}.
\end{equation*}%
Further, using $\mathcal{F}_{t_{k}}$-measurability of $\rho _{k}$ and the
conditional variants of (\ref{Ba06-lp}), (\ref{Ba10-lp}) and (\ref{Ba11-lp})
and the Cauchy-Bunyakovsky inequality (twice), we get for $p\geq 2:$%
\begin{eqnarray}
&&\mathbb{E}\left( \left\vert S_{k+1}\right\vert ^{2p-2}-\left\vert \rho
_{k}\right\vert ^{2p-2}\right) (\rho _{k},r_{k+1})  \label{thm42} \\
&\leq &K\mathbb{E}\left\vert \rho _{k}\right\vert ^{2p-1}{%
(1+|X(t_{k})|^{2\varkappa -2}+|X_{k}|^{2\varkappa -2})^{1/4}}%
h^{q_{2}+1/2}(1+|X_{k}|^{2\alpha })^{1/2}.  \notag
\end{eqnarray}

Due to $\mathcal{F}_{t_{k}}$-measurability of $\rho _{k},$ the conditional
variants of (\ref{Ba06-lp}) and (\ref{Ba11-lp}) and the Cauchy-Bunyakovsky
inequality (twice), we obtain for the third term on the right-hand side of (%
\ref{thm3}):%
\begin{gather}
\mathbb{E}\left\vert S_{k+1}\right\vert ^{2p-2}(Z_{k+1},r_{k+1})\leq \mathbb{%
E}[\mathbb{E}\left( S_{k+1}|^{4p-4}|\mathcal{F}_{t_{k}}\right) ^{1/2}\mathbb{%
E}\left( |Z_{k+1}|^{4}|\mathcal{F}_{t_{k}}\right) ^{1/4}\mathbb{E}\left(
|r_{k+1}|^{4}|\mathcal{F}_{t_{k}}\right) ^{1/4}]  \label{thm5} \\
\leq K\mathbb{E}\left\vert \rho _{k}\right\vert ^{2p-1}{(1+|X(t_{k})|^{2%
\varkappa -2}+|X_{k}|^{2\varkappa -2})^{1/4}}h^{q_{2}+1/2}(1+|X_{k}|^{4%
\alpha })^{1/4}.  \notag
\end{gather}

Due to $\mathcal{F}_{t_{k}}$-measurability of $\rho _{k}$ and due to the
conditional variants of (\ref{Ba06-lp}) and (\ref{Ba10-lp}) and the
Cauchy-Bunyakovsky inequality, we estimate the third term on the right-hand
side of (\ref{thm1}): 
\begin{eqnarray}
K\sum_{l=2}^{2p}\mathbb{E}\left\vert S_{k+1}\right\vert ^{2p-l}|r_{k+1}|^{l}
&\leq &K\sum_{l=2}^{2p}\mathbb{E}[\mathbb{E}(\left\vert S_{k+1}\right\vert
^{4p-2l}|\mathcal{F}_{t_{k}})^{1/2}\mathbb{E}(|r_{k+1}|^{2l}|\mathcal{F}%
_{t_{k}})^{1/2}]  \label{thm6} \\
&\leq &K\sum_{l=2}^{2p}\mathbb{E}[|\rho
_{k}|^{2p-l}h^{lq_{2}}(1+|X_{k}|^{2l\alpha })^{1/2}].  \notag
\end{eqnarray}

Substituting (\ref{thm2})-(\ref{thm6}) in (\ref{thm1}) and recalling that $%
q_{1}\geq q_{2}+1/2,$ we obtain%
\begin{gather*}
\mathbb{E}|\rho _{k+1}|^{2p}\leq \mathbb{E}|\rho _{k}|^{2p}(1+Kh)+K\mathbb{E}%
\left\vert \rho _{k}\right\vert ^{2p-1}(1+|X_{k}|^{2\alpha
})^{1/2}h^{q_{2}+1/2} \\
+K\mathbb{E}\left\vert \rho _{k}\right\vert ^{2p-1}{(1+|X(t_{k})|^{2%
\varkappa -2}+|X_{k}|^{2\varkappa -2})^{1/4}}h^{q_{2}+1/2}(1+|X_{k}|^{2%
\alpha })^{1/2} \\
+K\mathbb{E}\left\vert \rho _{k}\right\vert ^{2p-1}{(1+|X(t_{k})|^{2%
\varkappa -2}+|X_{k}|^{2\varkappa -2})^{1/4}}h^{q_{2}+1/2}(1+|X_{k}|^{4%
\alpha })^{1/4} \\
+K\sum_{l=2}^{2p}\mathbb{E}[|\rho _{k}|^{2p-l}h^{lq_{2}}(1+|X_{k}|^{2\alpha
l})^{1/2}] \\
\leq \mathbb{E}|\rho _{k}|^{2p}(1+Kh)+K\mathbb{E}\left\vert \rho
_{k}\right\vert ^{2p-1}{(1+|X(t_{k})|^{2\varkappa -2}+|X_{k}|^{2\varkappa
-2})^{1/4}}h^{q_{2}+1/2}(1+|X_{k}|^{2\alpha })^{1/2} \\
+K\sum_{l=2}^{2p}\mathbb{E}[|\rho _{k}|^{2p-l}h^{lq_{2}}(1+|X_{k}|^{2l\alpha
})^{1/2}].
\end{gather*}%
Then using Young's inequality and the conditions (\ref{Xmom}) and (\ref%
{appmomt}), we obtain 
\begin{equation*}
\mathbb{E}|\rho _{k+1}|^{2p}\leq \mathbb{E}|\rho _{k}|^{2p}+Kh\mathbb{E}%
|\rho _{k}|^{2p}+K{(1+\mathbb{E}|X}_{0}{|^{\beta p(\varkappa -1)+2p\alpha
\beta })}h^{2p(q_{2}-1/2)+1}
\end{equation*}%
whence (\ref{Ba08-lp}) with integer $p\geq 1$ follows by application of
Gronwall's inequality. Then by Jensen's inequality (\ref{Ba08-lp}) holds for
non-integer $p$ as well. $\square $

\section{Proof of Lemma~\protect\ref{lm:Bat02-lp} \label{sec:prooflem}}

Lemma~\ref{lm:Bat02-lp} is an analogue of Lemma~1.1.3 in \cite{MT6}.

\noindent \textbf{Proof}. Introduce the process $%
S_{t,x,y}(s)=S(s):=X_{t,x}(s)-X_{t,y}(s)$ and note that $Z(s)=S(s)-(x-y)$.
We first prove (\ref{Ba10-lp}). Using the Ito formula and the condition (\ref%
{olc2}) (recall that (\ref{olc2}) implies (\ref{Xmom})), we obtain for $%
\theta \geq 0:$%
\begin{gather*}
\mathbb{E}|S(t+\theta )|^{2p}=|x-y|^{2p}+2p\int_{t}^{t+\theta }\mathbb{E}%
|S|^{2p-2}%
\Bigg[%
S^{\intercal }(a(t,X_{t,x}(s))-a(t,X_{t,y}(s)))%
\Bigg.
\\
\left. +\frac{1}{2}\sum_{r=1}^{m}|\sigma _{r}(t,X_{t,x}(s))-\sigma
_{r}(t,X_{t,y}(s))|^{2}\right] \,ds \\
+2p(p-1)\int_{t}^{t+\theta }\mathbb{E}|S|^{2p-4}\left\vert S^{\intercal
}(s)\sum_{r=1}^{m}[\sigma _{r}(t,X_{t,x}(s))-\sigma
_{r}(t,X_{t,y}(s))]\right\vert ^{2}ds \\
\leq |x-y|^{2p}+2p\int_{t}^{t+\theta }\mathbb{E}|S|^{2p-2}%
\Bigg[%
S^{\intercal }(a(t,X_{t,x}(s))-a(t,X_{t,y}(s)))%
\Bigg.
\\
\left. +\frac{2p-1}{2}\int_{t}^{t+\theta }\sum_{r=1}^{m}|\sigma
_{r}(t,X_{t,x}(s))-\sigma _{r}(t,X_{t,y}(s))|^{2}\right] \,ds \\
\leq |x-y|^{2p}+2pc_{1}\int_{t}^{t+\theta }\mathbb{E}|S(s)|^{2p}\,ds
\end{gather*}%
from which (\ref{Ba10-lp}) follows after applying Gronwall's inequality.

Now we prove (\ref{Ba11-lp}). Using the Ito formula and the condition (\ref%
{olc2}), we obtain for $\theta \geq 0:$ 
\begin{gather}
\mathbb{E}\left\vert Z(t+\theta )\right\vert ^{2p}=2p\int_{t}^{t+\theta }%
\mathbb{E}|Z|^{2p-2}%
\Bigg[%
Z^{\intercal }(a(t,X_{t,x}(s))-a(t,X_{t,y}(s)))%
\Bigg.
\label{lmx1} \\
\left. +\frac{1}{2}\sum_{r=1}^{m}|\sigma _{r}(t,X_{t,x}(s))-\sigma
_{r}(t,X_{t,y}(s))|^{2}\right] \,ds  \notag \\
+2p(p-1)\int_{t}^{t+\theta }\mathbb{E}|Z|^{2p-4}\left\vert Z^{\intercal
}\sum_{r=1}^{m}[\sigma _{r}(t,X_{t,x}(s))-\sigma
_{r}(t,X_{t,y}(s))]\right\vert ^{2}\,ds  \notag
\end{gather}%
\begin{gather*}
\leq 2p\int_{t}^{t+\theta }\mathbb{E}|Z|^{2p-2}(s)%
\Bigg[%
S^{\intercal }(a(t,X_{t,x}(s))-a(t,X_{t,y}(s)))%
\Bigg.
\\
\left. +\frac{2p-1}{2}\int_{t}^{t+\theta }\sum_{r=1}^{m}|\sigma
_{r}(t,X_{t,x}(s))-\sigma _{r}(t,X_{t,y}(s))|^{2}\right] \,ds \\
-2p\int_{t}^{t+\theta }\mathbb{E}%
|Z|^{2p-2}(x-y,a(t,X_{t,x}(s))-a(t,X_{t,y}(s)))ds
\end{gather*}%
\begin{equation*}
\leq 2pc_{1}\int_{t}^{t+\theta }\mathbb{E}|Z|^{2p-2}|S|^{2}\,ds-2p%
\int_{t}^{t+\theta }\mathbb{E}%
|Z|^{2p-2}(x-y,a(t,X_{t,x}(s))-a(t,X_{t,y}(s)))ds.
\end{equation*}%
Using Young's inequality, we get for the first term in the right-hand side
of (\ref{lmx1}):%
\begin{gather}
2pc_{1}\int_{t}^{t+\theta }\mathbb{E}|Z|^{2p-2}|S|^{2}\,ds\leq
4pc_{1}\int_{t}^{t+\theta }\mathbb{E}|Z|^{2p-2}(|Z|^{2}+|x-y|^{2})\,ds
\label{lmx2} \\
\leq K\int_{t}^{t+\theta }\mathbb{E}|Z|^{2p}ds+K|x-y|^{2}\int_{t}^{t+\theta }%
\mathbb{E}|Z|^{2p-2}ds.  \notag
\end{gather}

Consider the second term in the right-hand side of (\ref{lmx1}). Using
Hoelder's inequality (twice), (\ref{olc3}), (\ref{Ba10-lp}) and (\ref{Xmom}%
), we obtain%
\begin{eqnarray}
&&-2p\int_{t}^{t+\theta }\mathbb{E}%
|Z|^{2p-2}(x-y,a(t,X_{t,x}(s))-a(t,X_{t,y}(s)))ds  \label{lmx3} \\
&\leq &2p\int_{t}^{t+\theta }\mathbb{E}%
|Z|^{2p-2}|a(t,X_{t,x}(s))-a(t,X_{t,y}(s))||x-y|ds  \notag \\
&\leq &K|x-y|\int_{t}^{t+\theta }\left[ \mathbb{E}|Z|^{2p}\right] ^{1-1/p}%
\left[ \mathbb{E}|a(t,X_{t,x}(s))-a(t,X_{t,y}(s))|^{p}\right] ^{1/p}ds 
\notag
\end{eqnarray}%
\begin{eqnarray*}
&\leq &K|x-y|\int_{t}^{t+\theta }\left[ \mathbb{E}|Z|^{2p}\right] ^{1-1/p} \\
&&\times (\mathbb{E}[(1+|X_{t,x}(s)|^{2\varkappa
-2}+|X_{t,y}(s)|^{2\varkappa
-2})^{p/2}|X_{t,x}(s)-X_{t,y}(s)|^{p}])^{1/p}\,ds \\
&\leq &K|x-y|\int_{t}^{t+\theta }\left[ \mathbb{E}|Z|^{2p}\right]
^{1-1/p}\left( \mathbb{E}[(1+|X_{t,x}(s)|^{2\varkappa
-2}+|X_{t,y}(s)|^{2\varkappa -2})^{p}]\right) ^{1/2p} \\
&&\times \left( \mathbb{E}[|X_{t,x}(s)-X_{t,y}(s)|^{2p}]\right) ^{1/2p}\,ds
\\
&\leq &K\left\vert x-y\right\vert ^{2}(1+|x|^{2\varkappa -2}+|y|^{2\varkappa
-2})^{1/2}\int_{t}^{t+\theta }\left[ \mathbb{E}|Z|^{2p}\right] ^{1-1/p}ds.
\end{eqnarray*}

Substituting (\ref{lmx2}) and (\ref{lmx3}) in (\ref{lmx1}) and applying
Hoelder's inequality to $\mathbb{E}|Z|^{2p-2}\cdot 1$, we get 
\begin{equation}
\mathbb{E}\left\vert Z(t+\theta )\right\vert ^{2p}\leq K\int_{t}^{t+\theta }%
\mathbb{E}|Z|^{2p}ds+K\left\vert x-y\right\vert ^{2}(1+|x|^{2\varkappa
-2}+|y|^{2\varkappa -2})^{1/2}\int_{t}^{t+\theta }\left[ \mathbb{E}|Z|^{2p}%
\right] ^{1-1/p}ds  \label{lmx4}
\end{equation}%
whence we obtain (\ref{Ba11-lp}) for integer $p\geq 1$ using Gronwall's
inequality as, e.g. in \cite[p. 360]{Gronw}, and then by Jensen's inequality
for non-integer $p>1$ as well. $\ \square $


\begin{thebibliography}{99}
\bibitem{Nawaf} N. Bou-Rabee, E. Vanden-Eijnden. A patch that imparts
unconditional stability to explicit integrators for Langevin-like equations. 
\textit{J. Comp. Phys.} \textbf{231} (2012), 2565--2580.

\bibitem{Gyo98} I. Gy\"{o}ngy. A note on Euler's approximations. \textit{%
Poten. Anal.} \textbf{8\ }(1998), 205--216.

\bibitem{Has-B80} R.Z. Hasminskii. \textit{Stochastic Stability of
Differential Equations}. Sijthoff \& Noordhoff, 1980.

\bibitem{HW} E. Hairer, G. Wanner. \textit{Solving Ordinary Differential
Equations II: Stiff and Differential-Algebraic Problems.} Springer, 2004.

\bibitem{HMS} D.J. Higham, X. Mao, A.M. Stuart. Strong convergence of
Euler-type methods for nonlinear stochastic differential equations. \textit{%
SIAM J. Num. Anal.} \textbf{40 }(2003), 1041--1063.

\bibitem{Hu} Y. Hu. Semi-implicit Euler-Maruyama scheme for stiff stochastic
equations. \textit{In:} \textit{Stochastic Analysis and Related Topics V}: 
\textit{The Silvri Workshop} (Ed.: H. Koerezlioglu), Progr. Probab. 38,
Birkhauser, Boston, 1996, pp. 183--202.

\bibitem{HutJen12} M. Hutzenthaler, A. Jentzen. \textit{Numerical
approximation of stochastic differential equations with non-globally
Lipschitz continuous coefficients}. Preprint, 2012.

\bibitem{HutJenKlo12} M. Hutzenthaler, A. Jentzen, P.E. Kloeden. Strong
convergence of an explicit numerical method for SDEs with non-globally
Lipschitz continuous coefficients. \textit{Ann. Appl. Probab.} \textbf{22}
(2012), 1611--1641.

\bibitem{KP} P.E. Kloeden, E. Platen. \textit{Numerical Solution of
Stochastic Differential Equations}. Springer, 1992.

\bibitem{Stua} J.C. Mattingly, A.M. Stuart, D.J. Higham. Ergodicity for SDEs
and approximations: Locally Lipschitz vector fields and degenerate noise. 
\textit{Stoch. Proc. Appl.} \textbf{101} (2002), 185--232.

\bibitem{8} G.N. Milstein. A theorem on the order of convergence of
mean-square approximations of solutions of systems of stochastic
differential equations. \textit{Theor. Prob. Appl.} \textbf{32} (1987),
738--741.

\bibitem{GN} G.N. Milstein. \textit{Numerical Integration of Stochastic
Differential Equations}. Kluwer Academic Publishers, 1995.

\bibitem{MPS98} G.N. Milstein, E. Platen, H. Schurz. Balanced implicit
methods for stiff stochastic systems. \textit{SIAM J. Num. Anal.} \textbf{35 
}(1998), 1010--1019.

\bibitem{MRT2} G.N. Milstein, Yu.M. Repin, M.V. Tretyakov. Numerical methods
for stochastic systems preserving symplectic structure. \textit{SIAM J. Num.
Anal.} \textbf{40 (}2002), 1583--1604.

\bibitem{MT6} G.N. Milstein, M.V. Tretyakov. \textit{Stochastic Numerics for
Mathematical Physics}. Springer, 2004.

\bibitem{GNT04} G.N. Milstein, M.V. Tretyakov. Numerical integration of
stochastic differential equations with nonglobally Lipschitz coefficients. 
\textit{SIAM J. Num. Anal.} \textbf{43\ (}2005), 1139--1154.

\bibitem{GNTerg} G.N. Milstein, M.V. Tretyakov. Computing ergodic limits for
Langevin equations. \textit{Physica D} \textbf{229} (2007), 81--95.

\bibitem{filter} G.N. Milstein, M.V. Tretyakov. Monte Carlo algorithms for
backward equations in nonlinear filtering. \textit{Adv. Appl. Prob.} \textbf{%
41} (2009), 63--100.

\bibitem{Gronw} D.S. Mitrinovic, J. Pecaric, A.M. Fink. \textit{Inequalities
Involving Functions and their Integrals and Derivatives}\emph{.} Kluwer,
1994.

\bibitem{Orteg} J. M. Ortega, W.C. Rheinboldt. \textit{Iterative Solution of
Nonlinear Equations in Several Variables}. SIAM, 2000.

\bibitem{SzpMao10} L.\ Szpruch, X.\ Mao. Strong convergence and stability of
numerical methods for non-linear stochastic differential equations under
monotone condition. \textit{J. Comp. App. Math.} \textbf{238} (2013), 14--28.

\bibitem{Tal99} D. Talay. Stochastic Hamiltonian systems: exponential
convergence to the invariant measure, and discretization by the implicit
Euler scheme. \textit{Markov Proc. Relat. Fields} \textbf{8} (2002),
163--198.
\end{thebibliography}
\end{document}